%&amstex
\input amstex.tex
\documentstyle{amsppt}
\magnification=1200
\hsize=150truemm
\vsize=224.4truemm
\hoffset=4.8truemm
\voffset=12truemm

\TagsOnRight
\NoBlackBoxes
\NoRunningHeads

\refstyle{A}

\topmatter

\title Harmonic cocycles and cohomology of arithmetic groups (in positive
characteristic).
\endtitle

\author Marc Reversat \endauthor

\affil Laboratoire de Math\'ematiques \'Emile Picard, Unit\'e Mixte de Recherche
Universit\'e-C.N.R.S. 5580. 
\endaffil

\address {Universit\'e Paul Sabatier, 118 route de Narbonne, 31062 Toulouse c\'edex 4. France}
\endaddress

\email reversat\@picard.ups-tlse.fr \endemail

\date{September 1999}
\enddate

\abstract{Let $K$ be a global field of characteristic $p>0$. We study the cohomology of
arithmetic subgroups $\Gamma $ of $SL_{n+1}(K)$ (with respect to a fixed place of $K$), under the
hypothesis that these groups have no $p'$-torsion (any arithmetic group possesses a normal
subgroup of finite index without $p'$-torsion). We define the cohomology of $\Gamma $
with compact supports and values in ${\Bbb Z}[1/p]$, and we relate it to spaces of harmonic
cocycles, also with compact supports (\S 3). We give a description of the locus of these
supports, in particular by introducing a notion of cusp in dimension $n\geq 1$ (\S 4) and we
calculate``geometrically" the Euler-Poincar\'e characteristic of this cohomology, up to torsion
(\S 5).} \endabstract

%\thanks...\endthanks ;  \keywords...\end.. ;  \subjclass...\end.. 

\endtopmatter

\document

\head 1 Introduction.\endhead

\noindent {\bf (1.0.1)}  
Let $K$ be a global field of characteristic $p>0$, i.e. the function field of a geometrically
irreducible curve defined over a finite field of characteristic $p$. One sets $G=SL_{n+1}$, where
$n\geq 1$ is an integer. Let $\infty$ be a fixed place of $K$, $K_{\infty}$ be the completion of
$K$ at this place and $\frak I$ be the the Bruhat-Tits building of $G(K_\infty )$ (see \S 2 for 
references and a brief description of this building); $\frak I$ is a contractible simplicial
complex of dimension $n$, where $G(K_\infty)$ acts simplicially. Let $\Gamma $ be an arithmetic
subgroup of
$G(K_\infty)$ (see (2.4.1) for a precise definition), the $G(K_\infty)$-action induces a 
$\Gamma $-action on $\frak I$. 

\vskip5pt
The purpose of this paper is to give information about the cohomology of $\Gamma $,
more precisely, to relate this cohomology to geometric or simplicial properties of the building
$\frak I$ or of its quotient by $\Gamma $.\par
  
\vskip5pt  Links of such kind have been studied particularly for discrete groups\break 
$\Gamma \subset G(K_\infty)$ which are {\it cocompact}, i.e. such that $\Gamma \backslash \frak I$ is
a finite complex, and, moreover, which are {\it torsion free}, and thus which {\it act freely} 
on $\frak I$. This has been carried out by A. Borel and J.-P. Serre in \cite {Bo-Se 2}, by G.A.
Mustafin in
\cite {Mus}, and, in their study of the Drinfeld symmetric spaces (which are rigid analytic
realizations of the simplicial complexes $\frak I$), by  P. Schneider and U. Stuhler (\cite
{Sc-St}). When the characteristic of the base field is zero, i.e. in the $p$-adic case or in the
classical case, many authors have studied the cohomology of arithmetic, torsion free and cocompact
groups acting on
$\frak I$ or on other symmetric spaces: see for instance the three papers mentioned above, 
in characteristic zero, and J.-P. Serre (\cite {Se 1}), H. Garland (\cite {GarH}), A. Borel
and J.-P. Serre (\cite {Bo-Se 1}) etc., and a recent paper of E. de Shalitt (\cite {dS 2}). There is
also, only in dimension $1$ but in all characteristics, the study of Schottky groups by Mumford
curves (see
\cite {Mum} and \cite {Ger-vdP}).\par
  In the situation studied here, the arithmetic groups are discrete, but {\it they are not
cocompact} (the complex $\Gamma \backslash {\frak I}$ is not finite) and {\it they have torsion
elements} (simplices of $\frak I$ have non trivial stabilizers in $\Gamma $). When $n=1$, i.e.
in the dimension $1$ case, this has been extensively studied, following the seminal work of
Drinfeld \cite {Dr} (see \cite {Se 2} ch. 2, \cite {Gek-Re} and \cite {A-B}). 
To explain what we do here, we now recall somes properties of the dimension $1$ case. 

\subhead 1.1 The dimension $1$ case \endsubhead 

Let $\Gamma $ be an arithmetic subgroup of
$G(K)=SL_2(K)$, where $K$ is a global field of characteristic $p>0$.

\vskip5pt \noindent {\bf (1.1.1)} In the dimension $1$ case, the building $\frak I$ is a tree and
the quotient $\Gamma \backslash {\frak I}$ is a connected graph which is the union of a finite
graph
$(\Gamma \backslash {\frak I})^\circ$ without ends and of finitely many half lines; the graph
$(\Gamma \backslash {\frak I})^\circ$ contains all the homology of $\Gamma \backslash {\frak I}$
and the half lines corresponds to the cusps of $\Gamma $ (see \cite {Se 2} ch. 2, \cite {vdP},
see also \cite {Gek-Re} \S 2).

\vskip5pt \noindent {\bf (1.1.2)} The harmonic cocycles are functions $f$ defined on the set of
oriented edges of $\frak I$, with for instance values in $\Bbb Z$, such 
that for any oriented edge $e$, $f(-e)=-f(e)$, and such that for any vertex $v$
of $\frak I$, $\sum f(e)=0$, where $e$ runs in the set of edges beginning at $v$. An important fact
is that the group ${\underline H}_!({\Bbb Z})^\Gamma $, 
of $\Bbb Z$-valued harmonic cocycles invariant under
$\Gamma $ and with compact (finite) supports modulo $\Gamma $, is canonically isomorphic with
$H^1(\Gamma ,{\Bbb Z})$, under the hypothesis that $\Gamma $ has no $p'$-torsion, i.e. that any
torsion element of $\Gamma $ is of order a power of $p$
(\cite {Gek-Re} prop. 3.4.5). This last hypothesis is indeed natural and not too restrictive,
because any arithmetic group contains a normal subgroup of finite index and having this property
(this is right in all dimensions, see 2.4).\par
  The group $H^1(\Gamma ,{\Bbb Z})$ can also be interpreted as the cohomology group of the cochain
complex $(C^i({\Bbb Z}))_{i=0,1}$ (definitions are given in 
\S 3.1, see \cite {Se 2} ch. 2, \S 2.8).

\vskip5pt \noindent {\bf (1.1.3)} The graph $(\Gamma \backslash {\frak I})^\circ$ has no ends;
let ${\overline T}$ be one of its maximal subtrees and $T$ be a connected subtree of $\frak I$ wich
is a lift of ${\overline T}$. Let $e_1,\cdots ,e_g$be edges of $\frak I$, beginning in $T$ and
equal modulo $\Gamma $ to the edges of $((\Gamma \backslash {\frak I})^\circ - {\overline T})$.
Let $\gamma_1,\cdots ,\gamma_g$ be elements of $\Gamma $ such that the end of $\gamma_i(e_i)$ is
in $T$ ($1\leq i\leq g$). Then, the group $G:=<\gamma_1,\cdots ,\gamma_g>$ is free , it is
canonically identified with the fundamental group of $(\Gamma \backslash {\frak I})^\circ$ and 
$$G\subset \Gamma \rightarrow \Gamma /\Gamma_{\text tors}$$ 
is an isomorphism, where $\Gamma_{\text tors}$ is the subgroup of $\Gamma $ generated by its
torsion elements (\cite {Re}, compare with \cite {Se 2} ch. 1, \S 3.3). As a consequence (but it can
be seen more directly, see \cite {Gek-Re} \S 3.2) one has for the Euler-Poincar\'e characteristic
of $\Gamma $ with values in $\Bbb Q$
$$\chi (\Gamma ,{\Bbb Q})=1-g.$$

\subhead 1.2 Summary of results \endsubhead

The object of this article is to extend to the case of dimension $n>1$ some of the properties
reviewed in the preceding paragraph.\par
  We suppose that the arithmetic group $\Gamma \subset G(K)$ has no $p'$-torsion.
Let $R$ be a (commutative) ring such that multiplication by $p$ is invertible. 
Let $C^\cdot (R)$ be the complex of $R$-valued cochains, i.e. 
$$C^q(R)={\text {Hom}}_{\Bbb Z}({\Bbb Z}[{\frak I}^\star_q], R) \ , \ 0\leq q\leq n,$$
where ${\frak I}^\star_q$ is the set of oriented $q$- simplices of ${\frak I}$. We first prove
that, canonically
$$H^\cdot (\Gamma ,R) \simeq H^\cdot (C^\cdot(R)^\Gamma )$$
(proposition (3.1.4)), indeed this is a direct application of an argument of J.P. Serre: \cite {Se
1}, p.95; here, the hypothesises that $\Gamma $ has no $p'$-torsion and that $p$ is invertible in 
${\text {End}}_{\Bbb Z}(R)$ are very useful.
This formula leads to the definition of {\it cohomology of $\Gamma $ with compact
supports}, which is the part in $H^\cdot (\Gamma ,R)$ corresponding in $H^\cdot (C^\cdot(R)^\Gamma
)$ to the cocycles having, up to coboundaries, finite supports modulo $\Gamma $ (definition
(3.1.5)). Afterwards, we introduce the notion of harmonic cocycle (definition (3.2.5)), as in the
paper 
\cite {Ga} of H. Garland
\footnote{H. Garland writes that the idea of this notion goes back to
Eckmann and Hodge. A variant, especially for the maximal rank (harmonic cocycles defined on pointed
chambers), is used in \cite {Sc-Te 1}, \cite {Sc-Te 2}, \cite {dS 1} and \cite {dS 2}.}. Let
$${\underline H}^q_!(R)^\Gamma  \ {\text {and}} \  H^q_! (\Gamma ,R) \ , \ 0\leq q\leq n$$
denote respectively 
the $R$-module of harmonic cocycles defined on the oriented $q$-simplices, invariant under $\Gamma$
 and with (finite)  compact supports modulo $\Gamma $, and the $q$-th cohomological $R$-module 
of $\Gamma $ with compact
supports. Set $R={\Bbb Z}[1/p]$ (indeed what we will do works for any subring of $\Bbb R$
containing ${\Bbb Z}[1/p]$).
We prove, among other things, that the module of harmonic cocycles
${\underline H}^q_!(R)^\Gamma$ is canonically isomorphic
to the torsion free submodule of $H^q_! (\Gamma ,R)$ (theorem (3.3.1)). This extends to the higher
dimensional case some aspects of the results recalled in (1.1.2); we will explain below (see
(1.3.1)) why we obtain an isomorphism with the cohomology of $\Gamma $ with compact supports,
instead of the whole cohomology, as in the dimension $1$ case.\par
  In the 4-th paragraph, we introduce a notion of cusp for $\Gamma $. One knows that 
$\Gamma \backslash {\frak I}$ is the ``skeleton" of a modular variety $M_\Gamma $ (\cite {Dr}, see
also \cite {De-Hu}), which is affine; from that point of view, the
cusps are skeletons of neighborhoods of the missing part, 
i.e. of the difference between $M_\Gamma $ and its canonical
completion (but we shall not say anything in this paper about modular varieties). It is not
difficult to prove that there are only finitely many cusps modulo $\Gamma $ (and that they
are rational over $K$, as in the dimension
$1$ case, see \cite {Se 2} ch. 2, \S 2.3). They are represented modulo $\Gamma $ by
finitely many sectors of $\frak I$, say ${\Cal S}_1,\cdots , {\Cal S}_d$, and the complementary
complex in $\frak I$ of 
$$\bigcup_{1\leq i\leq d}\Gamma {\Cal S}_i$$ 
is finite modulo $\Gamma $ (theorem
(4.2.2)). Thus, this complementary complex and the cusps, modulo $\Gamma $, give respectively an
analogy with $(\Gamma \backslash {\frak I})^\circ $ and the half lines (see (1.1.1)). The proof of
theorem (4.2.2) is long, the argument is based on results of G. Harder 
(s\"atze 2.1.1 and 2.2.2 of \cite {Ha 1}), concerning the
locus of supports of cusp forms (\cite {Ha 2} theorem 1.2.1).\par
  The 5-th paragraph gives a partial generalization of (1.1.3), in calculating ``geometrically"
the Euler-Poincar\'e characteristic of $\Gamma $, for the cohomology with values in $\Bbb Q$. This
geometric calculation is a step towards notions analogous to the group $G$ of (1.1.3); problems,
going in that way, are suggested in (5.0.5).

\subhead 1.3 Some comments \endsubhead

 \noindent {\bf (1.3.1)} Let $\Gamma $ be an arithmetic group whithout $p'$-torsion. We have obtained
a comparison theorem between the spaces 
${\underline H}_!^\cdot ({\Bbb Q})^\Gamma$ of harmonic cocycles invariant under $\Gamma $, with
finite supports modulo $\Gamma$ 
and with values in $\Bbb Q$ (we take here just $\Bbb Q$ to be short) and the spaces
$H_!^\cdot (\Gamma ,{\Bbb Q})$ of cohomology with compact supports (theorem (3.3.1)). In the one
dimensional case, $H_!^\cdot (\Gamma ,{\Bbb Q})$ is replaced by the ordinary cohomology
$H^\cdot (\Gamma ,{\Bbb Q})$. The reason for this difference is that, in the dimension $1$ case,
the ``simplicial complex at infinity" of  $\Gamma \backslash {\frak I}$  (which comes from the
spherical building of $\frak I$ at infinity, but which is {\it not} simply the quotient of this
spherical building by $\Gamma $),
chambers of which can be viewed as the cusps, is a finite complex of dimension $0$ (see
(1.1.1)), then there is no complexity (no cycle) at infinity. In that
situation, writing $P$ for the inverse image in $\frak I$ of the half lines of 
$\Gamma \backslash {\frak I}$ 
(the cusps, see (1.1.1)), a Mayer-Vietoris argument applied to the decomposition
$${\frak I}=({\frak I}-P)\cup P$$
yields equality between the ordinary cohomology of $\Gamma $ and 
the cohomology with compact supports.
When the dimension $n$  is greater than $1$, the simplicial complex 
at infinity of $\Gamma \backslash {\frak I}$ 
is still finite, but it has dimension $n-1>0$ and it may have non trivial homology; in this
case the two cohomologies of $\Gamma $ are not
equal. In our opinion, to give a description of the homological nature of the simplicial complex 
at infinity of 
$\Gamma \backslash {\frak I}$ is a very important question, which seems to underly the
conjecture (13.4.1) of \cite {Lau}.
It is related, as the next question, to calculations and interpretations of
Eisenstein cohomologies (see \cite {Ha 3} and op.cit.).

\vskip5pt \noindent {\bf (1.3.2)} In the dimension $1$ case, harmonic cocycles invariant under 
$\Gamma $ and with finite support modulo $\Gamma $ can be interpreted canonically,
for a suitable subgroup of the ad\`ele group of $G(K)$, as cusp forms that transform like the
special representation (\cite {Dr}, \cite {De-Hu}, \cite {Ge-Re}, \cite {A-B}, etc.). This is not so
clear in dimension greater than $1$. This problem is mentioned in the introduction of  \cite {dS 2}
and in \cite {Sc-St}, remark (c) p. 84. Our opinion is also that this is an
interesting question.

\vskip5pt \noindent {\bf (1.3.3)} Problems can be formulated to determine objects, in dimension
greater than $1$, similar to $G=<\gamma_1, \cdots , \gamma_g>$ (see (1.1.3)), they need
notations and definitions, which are introduced in \S 5, then we have written this problem
in (5.0.5).

\head 2 The Bruhat-Tits building of $SL_{n+1}(K_\infty)$ (review). Notations.\endhead

\subhead 2.1 Notations\endsubhead

  Let $K$ denote a global field of positive characteristic $p$, i.e. the function field
of a geometrically irreducible curve $\Cal C$ defined over a finite field $\Bbb F$ of
characteristic $p$; we suppose that $\Cal C$ is smooth and projective. Let $\infty$ be a
fixed place of $K$, i.e. a closed point of $\Cal C$. One sets
$$A=H^0({\Cal C}-\{\infty \}, {\Cal O}_{\Cal C}), \leqno {\text {\bf (2.1.1)}}$$
$A$ is the ring of elements of $K$ regular outside $\infty$.

\vskip5pt \noindent {\bf (2.1.2)} Let $\frak V$ be the set of places of $K$ (the set of
closed points of $\Cal C$) and let
${\frak V}_f={\frak V}-\{\infty \}$ the set of ``finite places" 
\footnote{This terminology comes from an analogy with the classical case: $A$ and $K$ look like
respectively $\Bbb Z$ and $\Bbb Q$, the place $\infty$ like the archimedean place of 
$\Bbb Q$...} of $K$.\par
 Let $\omega $ be a place of $K$, we denote by $K_\omega $ 
the completion of $K$ at $\omega $, by ${\Cal O}_\omega $ its valuation ring, by 
${\Bbb F}(\omega )$ the residue field and by $\pi_\omega$ a uniformizing parameter at 
$\omega$, we suppose that $\pi_\omega$ is in $K$.\par
 Usually, we will denote by the same letter a place and the normalized valuation of $K$
corresponding to it, except for $\infty$, this valuation being denoted $\omega_\infty$. If
$\omega$ is a place of $K$, the corrresponding normalized valuation satisfies 
$\omega (\pi_\omega )=[{\Bbb F}(\omega):{\Bbb F}]$.\par
  Let $\omega $ be a valuation of $K$, we denote by $\mid \ \mid_\omega $
(by $\mid \ \mid_\infty $ if $\omega =\infty $) the corresponding 
normalized absolute value, i.e. the absolute value which satisfies
$\mid \pi _\omega  \mid_\omega =\sharp ({\Bbb F}(\omega ))^{-1}$.

\vskip5pt \noindent {\bf (2.1.3)} Let $n\geq 1$ be an integer, one sets
$$G=SL_{n+1} \ ,$$
let $T$ be its standard torus (of diagonal matrices), $P$ and $U$ be respectively the standard
Borel subgroup (of upper triangular matrices) and its unipotent radical (of matrices with in
addition
$1$ on the diagonal). Let $N$ be the normalizer of $T$ in $G$. One sets\par
${\overline W}= N(K_\infty)/T(K_\infty)$ : the linear Weyl group ,\par
$W= N(K_\infty)/T({\Cal O}_\infty)$ : the affine Weyl group.\par
\noindent Let $B_\infty $ be the inverse image of $P({\Bbb F}(\infty ))$ by the natural map
$G({\Cal O}_\infty ) @>>> G({\Bbb F}(\infty ))$.\par
 Let $\Bbb A$ (resp. ${\Bbb A}_f$), be the ring of  ad\`eles (resp. finite ad\`eles) of
$K$, i.e. the set of elements
$(f_\omega )\in \prod K_\omega $, where $\omega $ runs in $\frak V$ (resp. ${\frak V}_f$), such
that $f_\omega \in {\Cal O}_\omega $ except for finitely many places $\omega $. 
The groups of ad\`eles and finite ad\`eles of $G(K)$ are denoted respectively $G({\Bbb A})$ and 
$G({\Bbb A}_f)$.

\vskip5pt \noindent {\bf (2.1.4)} {\it The groups of characters and of $1$-parameter
subgroups.} Let
$$ X^\star ={\text {Hom}}_{\text {gr}}(T(K _\infty ), K^\star _\infty ) \ \  {\text {and}}
 \ \ X_\star ={\text {Hom}}_{\text {gr}}(K^\star _\infty , T(K _\infty ))$$ 
be the groups of characters and of $1$-parameter subgroups of $T$ (we do not distinguish
between $T$, ${\Bbb G}_m$ and their $K_\infty$-valued points because these groups are split).
They are both free abelian groups of rank $n$. For any $i$, $0\leq i\leq n$, let 
$\chi_i \in X^\star$ be the character which maps 
$t={\text {diag}}(t_o,\cdots ,t_n)\in T(K^\star_\infty )$  to $t_i$, then
$$\{\chi_i \ / \ 0\leq i\leq n \}$$
generates $X^\star$ and any of its subsets with $n$ elements is a basis. For any $i$,\break
$0\leq i\leq n-1$, let 
$\lambda_i\in X_\star$ be the $1$-parameter subgroup which maps $z\in K^\star _\infty $
to $t={\text {diag}}(t_o,\cdots ,t_n)$ such that $t_j=1$ if $j\not= i,i+1$, $t_i=z$ and
$t_{i+1}=z^{-1}$. Then
$$\{\lambda_i \ / \ 0\leq i\leq n-1\}$$
is a basis of $X_\star$. There exists a perfect pairing
$$X^\star \ \times \ X_\star \ \rightarrow \ {\Bbb Z}$$
which to $(\chi , \lambda )$ associates the integer $m$ such that, for any $z$ in 
$K^\star _\infty $, $\chi (\lambda (z))=z^m$ (see for instance \cite {Bo} \S 8).\par
  Note that $N(K_\infty)$ acts on $X^\star$ and $X_\star$ via its action on $T(K_\infty)$ by
conjugation.

\subhead 2.2 The building $\frak I$ \endsubhead

  Now we recall some basic facts about the Bruhat-Tits building $\frak I$ of $G(K_\infty)$.
For generalities on complexes and buildings, see \cite {Ro}, \cite {Bro} and \cite {GarP}. The
building of
$G(K_\infty)$ is completely described in  \cite {Bru-Ti}, there is a good introduction to this
topic in \cite {Lan}.\par
   Let ${\Cal S}$ be a simplicial complex. We denote by ${\Cal S}_q$ the set of its simplices of
dimension $q$. Let ${\Cal S}'$ be a subcomplex of ${\Cal S}$; recall that, by definifion, 
${\Cal S}'$ is a union of closed simplexes,
 recall also that the complementary complex of ${\Cal S}'$ in ${\Cal S}$ is the minimal
subcomplex containing the complementary set of ${\Cal S}'$ in ${\Cal S}$, then,
maybe, ${\Cal S}'$ and its complementary complex have common simplices.

\vskip5pt \noindent {\bf (2.2.1)} {\it The fundamental apartment.} Let
$$V_0=X_\star \otimes_{\Bbb Z} {\Bbb R},$$
it is a real vector space with basis 
$$\{e_i=\lambda_i \otimes 1 \ / \ 0\leq i\leq n-1\}.$$
The pairing between $X^\star$ and $X_\star$, introduced in (2.1.4), induces a pairing between
$V_0$ and $X^\star \otimes_{\Bbb Z} {\Bbb R}$, denoted by $<\cdot ,\cdot>$; it gives 
 the identification
$$V_0^\star=X^\star \otimes_{\Bbb Z} {\Bbb R}$$
where $V_0^\star$ is the dual space of $V_0$. We set 
$$a_i=\chi_i \otimes 1 \ ,\ 0\leq i\leq n \ {\text {and}} \ 
a_{i,j}=a_i-a_j \ , \ 0\leq i,j\leq n,$$
the $a_{i,j}$'s, $i\not= j$, are the roots of $G(K_\infty)$ relative to $T(K_\infty)$ (\cite
{Bo} \S 8.17).\par
  The vector space $V_0$ is a Coxeter complex, the walls of the chambers being the hyperplanes
$$\{<a_{i,j}, \cdot >+k=0\} \ , \ 0\leq i\leq n \ , \  k\in {\Bbb Z}.$$
    
The torus $T(K_\infty)$ acts on $V_0$  by translations: the translation corresponding to 
$t={\text {diag}}(t_0,\cdots ,t_n)$ is given by the vector
$$\sum_{0\leq i\leq n-1}m_i e_i \ {\text {where}} \ 
m_i=-(\omega_\infty (t_0)+\cdots +\omega_\infty (t_i)) \ ,$$
then, for any $i$, $0\leq i \leq n-1$, and for any $z\in K^\star _\infty $ , the translation 
in $V_0$ corresponding to $\lambda_i(z)$ is given by the vector $-\omega _\infty (z)e_i$.  This
action gives a map
$T(K_\infty)/T({\Cal O}_\infty) @>>> V_0$ (because $T({\Cal O}_\infty)$ acts trivially on
$V_0$).\par
  The normalizer $N(K_\infty)$ acts on $V_0$, this action comes from the action on $X_\star$,
for this action $T(K_\infty)$ acts trivially on $V_0$. It follows a map\break
$N(K_\infty)/T(K_\infty) @>>> GL(V_0)$.\par
  Everything that we have said about the different actions on $V_0$ can be summarized in the
following diagram (compare with \cite {Lan}, proof of lemma 1.7)
$$
\CD
0 @>>> {{\frac {T(K_\infty)}{T({\Cal O}_\infty)}}\simeq {\Bbb Z}^n} @>>> 
                              {{\frac {N(K_\infty)}{T({\Cal O}_\infty)}}=W} @>>> 
                                    {{\frac {N(K_\infty)}{T(K_\infty)}}={\overline W}} @>>> 0  \\
@.      @VVV                  @VVV                     @VVV                                 @. \\
0 @>>>  V_0    @>>>       GL(V_0)\propto V_0  @>>>     GL(V_0)          @>>>                0
\endCD
$$
where the two arrows are exact, the second vertical map is induced by the two others, which are
explained above; $GL(V_0)\propto V_0$ is the affine transformation group of $V_0$.\par
  The space $V_0$ equipped with its simplicial structure and with the previous affine action of
$W$ (or of $N(K_\infty)$)
is called {\it the fundamental apartment of the building} $\frak I$; one now defines this building.

\vskip5pt \noindent {\bf (2.2.2)} {\it Definition of the building.} Let $i,j$ be integers such that
$0\leq i\not= j\leq n$ and let $U_{i,j}$ 
be the subgroups of $u=(u_{k,l})\in U(K_\infty)$ (see (2.1.3))
such that $u_{k,l}=0$ if $k\not= l$ and $(k,l)\not= (i,j)$. For any $x\in V_0$ let
$$U_{i,j,x}=\{u=(u_{k,l})\in U_{i,j} \ / \ \omega _\infty (u_{i,j})+<a_{i,j},x> \geq 0\}$$
and let $U_x$ be the group generated by the $U_{i,j,x}$'s, $0\leq i\not= j\leq n$.\par
  On $G(K_\infty)\times V_0$ one has the following equivalence relation (\cite {Bru-Ti} 2.5
and 7.4, see also \cite {Lan} \S 9): let $(g, x)$ and $(h, y)$ be in  $G(K_\infty)\times V_0$,
they are said to be equivalent (and one writes $(g, x)\sim (h, y)$) if there exists 
$n\in N(K_\infty)$ such that $y=nx$ and $g^{-1}hn\in U_x$. Let
$${\frak I} \ = \ {\frac {G(K_\infty)\times V_0}{\sim }}.$$ 
One has $V_0\subset {\frak I}$ canonically, by $x\mapsto {\text {Class}}(1, x)$. The group 
$G(K_\infty)$ acts on $\frak I$ (if $g,h\in G(K_\infty)$ and $x\in V_0$, the action of $g$ on 
${\text {Class}}(h, x)$ gives ${\text {Class}}(gh, x)$) and this action extends the one of $W$
on $V_0$. For any $x\in V_0\subset {\frak I}$, the stabilizer of $x$ in $G(K_\infty)$ is the
group
$$P_x:=\langle U_x \ , \ N_x\rangle \ {\text {where}} \ N_x=\{n\in N(K_\infty)  \ / \ nx=x\},$$
This is a bounded subgroup of $G(K_\infty)$, for the $\infty$-topology.\par 
  $\frak I$ has a simplicial structure, which is induced by that of $V_0$, indeed $\frak I$
is an affine building, its apartments are of the form $gV_0$, for $g\in G(K_\infty )$,
the action of $G(K_\infty)$ on $\frak I$ is simplicial. The building $\frak I$
is the {\it Bruhat-Tits building of $G(K_\infty )$}.

\vskip5pt \noindent {\bf (2.2.3)} {\it Another point of view.} We now recall the description of
the building $\frak I$ by ${\Cal O}_\infty$-submodules of $K^{n+1}_\infty$ (see \cite {GarP}
ch. 19).\par
  We mentioned above that $V_0$ is a Coxeter complex (see (2.2.1)). It follows from the
definition of the walls of the chambers that a point 
$v=\sum_{0\leq i\leq n-1} m_ie_i$ (with $m_i\in {\Bbb R}$) of $V_0$ is a vertex if and only
if for each $i,j$, $0\leq i,j\leq n$, $<a_{i,j}, v>\in {\Bbb Z}$.\par
  On $K^{n+1}_\infty$ one considers the following $G(K_\infty )$-action: for any
$z\in K^{n+1}_\infty$ and any $g\in G(K_\infty )$, the action of $g$ on $z$ is given by the
product of matrices $zg^{-1}$, $z$ being viewed as a one-line matrix. Note that, in \cite {GarP}
ch. 19, $gz$ is chosen, with $z$ viewed as a one-column matrix. Let 
$\{u_0,\cdots ,u_n\}$ be the canonical basis of $K^{n+1}_\infty$.\par
  Let $v=\sum_{0\leq i\leq n-1} m_ie_i$ be a vertex of $V_0$. We set 
$$M_v= \bigoplus _{0\leq i\leq n}\pi^{<a_{i,n}, v>}_\infty {\Cal O}_\infty u_i.$$
This is an ${\Cal O}_\infty$-submodule of $K^{n+1}_\infty$ of maximal rank $n+1$. Note that
for any $t={\text {diag}}(t_0,\cdots ,t_n)\in T(K_\infty )$, one has $M_{tv}=t(M_v)$, where,
on the left, $tv$ is given by the action of $t$ on $V_0$ (see (2.2.1)) and, on the right,
by the preceding action on  $K^{n+1}_\infty$.\par
  Let ${\frak I}'$ be the building of homothety classes $[M]$
of ${\Cal O}_\infty$-submodule $M$ of $K^{n+1}_\infty$ with rank $n+1$ (see \cite {GarP} ch.
19). Then, the map $v\mapsto [M_v]$ induces a simplicial isomorphism between  
$V_0$ in ${\frak I}$ and the apartement of ${\frak I}'$ with frame 
$\{{\Cal O}_\infty u_0,\cdots ,{\Cal O}_\infty u_n\}$ (op. cit.). With the
$G(K_\infty )$-action described above, it follows that there is a $G(K_\infty )$-isomorphism 
${\frak I}\simeq {\frak I}'$. 

\vskip5pt \noindent {\bf (2.2.4)} {\it The fundamental chamber.} 
  Let $C_0$ be the set of $x\in V_0$ such that
$$<a_{i,i+1}, x>>0 \ , \ 0\leq i\leq n-1 \ \ {\text {and}} \ \ <a_{n, 0}, x>+1>0.$$ 
By the isomorphism of (2.2.3), it corresponds in ${\frak I}'$ to the chamber having vertices
$[M_0],\cdots ,[M_n]$ with
$$
\align M_0&={\Cal O}_\infty u_0\oplus \cdots \oplus {\Cal O}_\infty u_n\\
       M_{i+1}&=\pi _\infty {\Cal O}_\infty u_0\oplus \cdots \oplus 
           \pi _\infty {\Cal O}_\infty u_i \oplus {\Cal O}_\infty u_i \oplus \cdots
             {\Cal O}_\infty u_n \ , \ 0\leq i\leq n-1.\\
\endalign
$$ 
More precisely, $C_0$ is a chamber of $\frak I$ having vertices $v_i$, 
$0\leq i\leq n$, such that $v_0$ is the origin $0$ of $V_0$ and $v_{i+1}$ is defined by
$$<a_{k,k+1}, v_{i+1}>=0 \ {\text {for all}} \ k, \ k\not= i, \ 0\leq k\leq n-1 \ {\text {and}} \ 
<a_{i,i+1}, v_{i+1}>=1.$$
One has $M_{v_i}=M_i$, $0\leq i\leq n$. 
It follows from \cite {GarP} ch. 19, \S 4, that the stabilizer of
$C_0$ in $G(K_\infty )$ is the parahoric subgroup $B_\infty$ (see (2.1.3)).\par
  This chamber $C_0$ is called {\it the fundamental chamber}, its vertex corresponding to the
origin $0\in V_0$ is called {\it the fundamental vertex}.

\subhead 2.3 The spherical building at infinity \endsubhead

  Let ${\frak I}_\infty $ be the spherical building at infinity of $\frak I$. 
We shall say almost nothing
about this building, for a complete description, see for instance  \cite {Bro} ch. VI, \S 9. As
generality, one just recalls that ${\frak I}_\infty $ is constructed by attaching to each
apartment of ${\frak I}$ a sphere at infinity; a point of ${\frak I}_\infty $ is an equivalent
class of half lines, for the relation of parallelism.\par
  The group $G(K_\infty)$ acts on ${\frak I}_\infty $ and this action is compatible with that
on ${\frak I}$. Let $\overline W$ be the linear Weyl group (see (2.1.3)), one knows that 
$\overline W$ is a finite reflexion group acting on $V_0$, this last one carrying the structure 
of Coxeter complex associated to $\overline W$ (\cite {Bro} ch. I and II).

\vskip5pt \noindent {\bf (2.3.1)} Let $v$ be a vertex of ${\frak I}$, $v\in V_0$, a {\it sector} of
${\frak I}$ in $V_0$, with vertex $v$, or beginning in $v$,
is the (topological) closure of $v+{\Cal C}$, where 
${\Cal C}\i V_0$ is a chamber for the simplicial stucture of $V_0$ coming from $\overline W$.
A sector of ${\frak I}$ is the image under the $G(K_\infty )$-action of a sector of ${\frak I}$
in $V_0$.

\vskip5pt \noindent {\bf (2.3.2)} {\it Remark.} {\bf (1)} A sector of ${\frak I}$ in $V_0$ is
a closed simplicial cone of $V_0$ (viewed as an apartment of ${\frak I}$). 
There is the following geometric description: 
take a chamber $C$ of ${\frak I}$ in $V_0$, let $v$ be a vertex of
$C$ and $H_1,\cdots ,H_n$ be the walls of $C$ containing $v$. Let $S$ be the closure of the 
connected component of
$$V_0 - H_1\cup \cdots \cup H_n$$
containing $C$. Then $S$ is a sector of ${\frak I}$ in $V_0$, with vertex $v$, and all these
sectors have this form.\par
 {\bf (2)} Our definition implies that a sector of ${\frak I}$ is a subcomplex of ${\frak I}$;
this is the reason why we insist that the sectors be closed, which
is not always the case in the literature (as in \cite {Bro}).

\vskip5pt \noindent {\bf (2.3.3)} Let $v$ be a vertex of ${\frak I}$ and let $\Sigma _v$ be the
set of sectors of ${\frak I}$ with vertex $v$. Then, there is a one to one map between 
$\Sigma _v$ and the set of closed chambers of ${\frak I}_\infty$, which to each $S\in \Sigma _v$
associates the closed chamber ${\overline \Delta }_\infty $ of ${\frak I}_\infty $ equal to the
set of all equivalent classes of half lines beginning in $v$ and contained in $S$ (\cite {Bro} ch.
VI, \S 9, lemma 2). One says that the closed chamber ${\overline \Delta }_\infty $ is defined by
the sector $S$ or that the sector $S$ ends at ${\overline \Delta }_\infty $; 
if $\Delta _\infty$ is the chamber of ${\frak I}_\infty $ with closure 
${\overline \Delta }_\infty $, one also says that $\Delta _\infty $ is defined by $S$ or that
$S$ ends at $\Delta _\infty $.

\vskip5pt \noindent {\bf (2.3.4)} Let $S$ and $S'$ be two sectors of ${\frak I}$ (with different
vertices) defining the same chamber $\Delta _\infty $ of ${\frak I}_\infty $, then there exists a
sector $S''$ contained in $S$ and $S'$ and defining the same chamber $\Delta _\infty $ of 
${\frak I}_\infty $ (op.cit. lemma 4).

\vskip5pt \noindent {\bf (2.3.5)} {\it The fundamental sector.} Let $H_i$ be the hyperplane of 
$V_0$, defined by $<a_{i,i+1}, \cdot >=0$, $0\leq i\leq n-1$. Let $S_0$ be the sector of $V_0$
which contains the fundamental chamber and which is the closure of a connected component of
$$V_0 - H_0\cup \cdots \cup H_{n-1}.$$
Then, the vertex of $S_0$ is the fundamental vertex $v_0$ and one has
$$S_0=\{x\in V_0 \ / \  <a_{i,i+1}, x>\geq 0 \ , \ 0\leq i\leq n-1\}.$$
This sector $S_0$ is called {\it the fundamental sector.}\par
  Let $C_\infty $ be the chamber of ${\frak I}_\infty $ defined by the fundamental sector $S_0$.
One calls this {\it the fundamental chamber at infinity}. 

\vskip5pt \noindent {\bf (2.3.6}
Let $K_\infty ^{n+1}$ be equipped by the $G(K_\infty)$-action as before (see (2.2.3)). There exists
a $G(K_\infty)$-isomorphism between the building ${\frak I}_\infty$
and the building of flags of $K_\infty ^{n+1}$, which is compatible with that
described in (2.2.3) (arguments to prove that are closed to those of (2.2.3), see
for instance \cite {Bro} ch.VI, \S 9). 
Via this isomorphism, the fundamental chamber $C_\infty$ corresponds to the flag with
vertices
$$K_\infty u_{n-i}\oplus \cdots \oplus K_\infty u_n \ \ , \ \ 0\leq i\leq n$$
($\{u_0,\cdots ,u_n \}$ is the canonical basis of $K_\infty ^{n+1}$). Then it is clear that the
stabilizer of $C_\infty $ in $G(K_\infty)$ is $P(K_\infty)$. There is a canonical 
$G(K_\infty)$-equivariant bijection of the set of chambers of ${\frak I}_\infty$ 
with $G(K_\infty)/P(K_\infty)$.

\subhead 2.4 Arithmetic groups \endsubhead

  A subgroup of $G(K)$ is called arithmetic if it is
commensurable with $G(A)$ (see (2.1.1)). Such a group $\Gamma $ is discrete in $G(K_\infty )$ for
the $\infty$-topology and the stabilizer in $\Gamma $ of any simplex of $\frak I$ is finite.\par
  For any arithmetic subgroup $\Gamma $ of $G(K)$, there exists an ideal $I$ of $A$,
$I\not= 0$, such that
$$ \Gamma (I):={\text {ker}}(G(A) @>{\text {canonical}}>> G(A/I))$$
is a subgroup of $\Gamma $ with finite index. The groups of the form $\Gamma (I)$ are whithout
$p$'-torsion, it means that any torsion element of $\Gamma (I)$ is of order a power of the
characteristic $p$ of the base field.

\head 3 Harmonic cocycles and cohomology. \endhead

  In all this paragraph, the letter $\Gamma $ denotes a subgroup of $G(K)$ which is arithmetic.
We suppose that\par

\vskip5pt \noindent {\bf (3.0.1)} {\it $\Gamma $ has no $p'$-torsion, i.e. that the order of any
torsion element of $\Gamma $ is a power of $p$, where $p$ is the characteristic of $K$.}

\vskip5pt Note that this hypothesis is not too restrictive, see 2.4.

\subhead 3.1 Cohomology of $\Gamma $\endsubhead

 Let $0\leq q\leq n$ be an integer
and let ${\frak I}^*_q$ be the set of {\it oriented} $q$-simplices of ${\frak I}$; recall that
a $q$-simplex is oriented if it is equipped with an equivalent class of total orderings of
its vertices, two orderings being equivalent if they differ by an even permutation of the
vertices. Let ${\Bbb Z}[{\frak I}^*_q]$ be the
 $\Bbb Z$-module generated by the {\it oriented} $q$-simplices of ${\frak I}$,
 with the relations, when $q\not= 0$,  $\sigma_1+\sigma_2=0$, if  $\sigma_1$ and $\sigma_2$
are the two  oriented $q$-simplices corresponding to the same non oriented $q$-simplex.
  For $1\leq q\leq n$, let 
$\partial_q: {\Bbb Z}[{\frak I}^*_q]\rightarrow {\Bbb Z}[{\frak I}^*_{q-1}]$ be the morphism
of\break
$\Bbb Z$-modules  which  maps the oriented simplex $[v_0,\cdots , v_q]$ ($v_0,\cdots , v_q$ are
its vertices) to
$\sum_{0\leq i\leq q}(-1)^i[v_0,\cdots ,\hat v_i,\cdots , v_q]$, where the hat in $\hat v_i$ means
that this term is omitted. Let $\varepsilon : {\frak I}_0(M)@>>>M$ be the augmentation map
(which to each finite sum $\sum_\sigma a_\sigma \sigma $, $\sigma \in {\frak I}_0(M)$ and
$a_\sigma \in {\Bbb Z}$, associates $\sum_\sigma a_\sigma $).\par
  Note that $\Gamma $ acts on each ${\Bbb Z}[{\frak I}^*_q]$ (by the restriction of the action
of $G(K_\infty)$) and that the maps $\partial_q$ are morphisms of ${\Bbb Z}[\Gamma ]$-modules.\par
  It is wellknown that the following sequence of ${\Bbb Z}[\Gamma ]$-modules

$$0@>>>{\Bbb Z}[{\frak I}^*_n]@>\partial_n >>\cdots @>\partial_{q+1} >>
{\Bbb Z}[{\frak I}^*_q]@>\partial_q >>
\cdots @>\partial_1 >>{\Bbb Z}[{\frak I}^*_0]@>\varepsilon >>{\Bbb Z}@>>>0
\leqno {\text {\bf (3.1.1)}}$$
is exact (because  ${\frak I}$ is contractible) and that it is split.

\vskip5pt \noindent {\bf (3.1.2)} Let $R$ be a $\Bbb Z$-module. We suppose that
{\it multiplication by the characteristic $p$ of $K$ is invertible in $End_{\Bbb Z}(R)$}. 
This hypothesis will be necessary in (3.1.4). For $0\leq q\leq n$, we set 
$C^q(R)=Hom_{\Bbb Z}({\Bbb Z}[{\frak I}^*_q], R)$; we suppose that $\Gamma $ acts trivially
on $R$, then  $C^q(R)$ is equiped with the usual action of $\Gamma $ 
$f\overset \gamma\to \mapsto f\circ \gamma^{-1}$ ($f$ in $C^q(R)$, $\gamma $ in $\Gamma $).  
Let $d^q: C^q(R)@>>>C^{q+1}(R)$ be the map defined by $f\mapsto f\circ \partial_{q+1}$, 
where $f$ belongs to $C^q(R)$ and
$0\leq q\leq n-1$; let $e : R@>>>C^0(R)$ be the map $r\mapsto r\varepsilon $, $r\in R$.  
They are maps of $R[\Gamma]$-modules.

\vskip5pt \noindent {\bf (3.1.3)} {\it The following sequence of $R[\Gamma]$-modules is exact}
$$0@>>>R@>e>>C^0(R)@>d^0 >>\cdots @>d^{q-1} >>C^q(R)@>d^q >>\cdots @>d^{n-1} >>C^n(R)@>>>0.$$

\vskip5pt This is a direct consequence of the fact that (3.1.1) is exact and split.

\proclaim{(3.1.4) Proposition} One has natural isomorphisms of $\Bbb Z$-modules
$$H^\cdot (\Gamma ,R)\simeq H^\cdot (C^\cdot (R)^\Gamma ).$$
\endproclaim
 The left hand side term is the cohomology of $\Gamma $ with value in the trivial\break 
$\Gamma $-module $R$, on the right this is the cohomology of the complex 
$C^\cdot (R)^\Gamma $ ($\{ \ \}^\Gamma $ is the subset of  $\{ \ \}$ of its elements invariant
under $\Gamma $).

\demo{Proof} The argument of the proof is from J.P. Serre (\cite {Se 1}, p.95). Consider an
injective resolution of $C^\cdot (R)$ in the category of complexes of $R[\Gamma ]$-modules
and apply the functor $H^0(\Gamma , \ )$. Consider
the usual two spectral sequences which tend both to the cohomology of the total complex.
It follows from (3.1.3) that the second spectral sequence degenerates in the $E_2$ term, which
implies that the cohomology of the total complex is $H^\cdot (\Gamma , R)$.\par
   On the other hand, we have 
$$C^q(R)= Hom_{\Bbb Z}({\Bbb Z}[{\frak I}^*_q], R)
=Hom_{\Bbb Z}(\coprod_\sigma {\Bbb Z}[\Gamma /\Gamma_\sigma ], R)
=\prod_\sigma Hom_{\Bbb Z}({\Bbb Z}[\Gamma /\Gamma_\sigma ], R)$$
where $\sigma $ runs in a representative set of ${\frak I}^*_q$ modulo $\Gamma $,
$\Gamma_\sigma $ is the stabilizer of $\sigma $ in $\Gamma $ and where
${\Bbb Z}[\Gamma /\Gamma_\sigma ]$ is the free ${\Bbb Z}$-module with basis 
$\Gamma /\Gamma_\sigma$. We also have
$$Hom_{\Bbb Z}({\Bbb Z}[\Gamma /\Gamma_\sigma ], R)\simeq
Hom_{{\Bbb Z}[\Gamma_\sigma ]} ({\Bbb Z}[\Gamma ], R)=coind^\Gamma _{\Gamma_\sigma }(R)$$
this isomorphism coming from the fact that the action of $\Gamma_\sigma $ on $R$ is trivial.
Then
$$H^\cdot (\Gamma , C^q(R))\simeq 
\prod_\sigma  H^\cdot (\Gamma , coind^\Gamma _{\Gamma_\sigma }(R))\simeq
\prod_\sigma  H^\cdot (\Gamma_\sigma , R)$$
the last isomorphism being given by the Shapiro's lemma. As $\Gamma $ is without $p'$-torsion
and because multiplication by $p$ is invertible in $End_{\Bbb Z}(R)$, 
one has $H^q (\Gamma_\sigma , R)=0$ for all $q>0$. It
follows that the first spectral sequence degenerates also in the $E_2$ term and then that the
cohomology of the total complex is that of $C^\cdot (R)^\Gamma $. \qed \enddemo

  Then, we are interested by the cohomology of the complex $C^q(R)^\Gamma $. The following
definition precises an important part of this cohomology.

\proclaim{(3.1.5) Definition} Let $R$ be a $\Bbb Z$-module such that multiplication by $p$ is
invertible in $End_{\Bbb Z}(R)$. Let $\Gamma $ be an arithmetic group. For any integer $q$, 
$0\leq q\leq n$, we denote by $Z_!^q(R)^\Gamma $ the set of
elements of $Z^q(R)^\Gamma = (Ker d^q)^\Gamma $ (see (3.1.3))  
which have finite supports modulo $\Gamma $; let $H_!^q(\Gamma ,R)$ be its image in 
$H^q(\Gamma ,R)$ by
$$Z_!^q(R)^\Gamma \subset Z^q(R)^\Gamma @>{\text {can. surj.}}>>H^q(C^\cdot (R)^\Gamma )
\simeq H^q(\Gamma ,R)$$
(see (3.1.4)).We say that $H_!^q(\Gamma ,R)$ is the $q$-th
$R$-module of cohomology with compact (finite) supports of the arithmetic group $\Gamma $.
\endproclaim

\vskip5pt \noindent {\bf (3.1.6) Remark.} The cohomology groups $H^{\cdot}(\Gamma ,R)$ are finitely
generated\break  $R$- modules. 
This is wellknown when $n=1$ (see e.g. \cite {Ge-Re}), when $n>1$, it comes from the fact
that the arithmetic groups are of finite type (\cite {Be 1}, see also \cite {Be 2}).

\subhead 3.2 Harmonic cocycles \endsubhead

\noindent {\bf (3.2.1)}  Let $s$ and $\sigma $ be two oriented simplices of ${\frak I}$.\par 
One says that $s$ is an oriented face of $\sigma $, and one writes
$s<\sigma $, if, as non oriented simplices, $s$ is a face of $\sigma $ and if the orientation of
$s$ is the restriction of that of $\sigma $.\par
 Write $\sigma=[v_0,\cdots , v_q]$, where $v_0,\cdots ,v_q$ are the vertices of $\sigma $ with a
numbering which represents their equivalent class of orderings. Then, an oriented face $s$ of
$\sigma$ of codimension one can be written $s=[v_0,\cdots ,\hat v_i,\cdots , v_q]$. We set
$\eta(s, \sigma)=(-1)^i$, it does not depend on the choice in an equivalent class
of the ordering of the vertices of $\sigma $.\par
       As before, $R$ is a $\Bbb Z$-module.\par

\proclaim{(3.2.2) Definition}  Let $g$ be an element of $C^q(R)$, with $1\leq q\leq n$, 
let $s\in {\frak I}_{q-1}^*$ be an oriented $q-1$-simplex. We set 
$$\delta^q(g)(s)=\sum_{\sigma\in {\frak I}_q^*, \sigma>s} \eta(s, \sigma)g(\sigma)$$
the function $\delta^q(g)$ is an element of $C^{q-1}(R)$ and $\delta^q$ is a 
$R$-linear and $G(K_\infty)$-linear map from $C^q(R)$ to $C^{q-1}(R)$, $1\leq q\leq n$.\par 
  We also set $\delta^0=0: C^0(R)@>>>R$. 
\endproclaim

\vskip5pt \noindent {\bf (3.2.3) Remark. (i)}   
Let $g$ be an element of $C^q(R)$, with $1\leq q\leq n$, let 
$s\in {\frak I}_{q-1}^*$ be an oriented $(q-1)$-simplex. Let $\sigma\in {\frak I}_q^*$ and suppose
that $s$ is an oriented face of $\sigma $. If we write $s=[v_0,\cdots , v_{q-1}]$, we have, up
to equivalence of orderings of the vertices, two possibilities for $\sigma $, the vertex $v$ of
$\sigma $ not in $s$ being fixed, namely
$\sigma=\sigma_1=[v, v_0,\cdots , v_{q-1}]$ and $\sigma=\sigma_2=[v_0, v,\cdots , v_{q-1}]$. We
have $\eta(s, \sigma_1)g(\sigma_1)=\eta(s, \sigma_2)g(\sigma_2)$. Then we can write
$$1/2(\eta(s, \sigma_1)g(\sigma_1)+\eta(s, \sigma_2)g(\sigma_2))=
\eta(s, \sigma_1)g(\sigma_1)=\eta(s, \sigma_2)g(\sigma_2)$$ 
and
$$1/2\sum_{\sigma\in {\frak I}_q^*, \sigma>s} \eta(s, \sigma)g(\sigma)$$  
makes sense (see the previous definition), but, with this formula as definition of $\delta^q$, the
following remark is wrong.\par
  {\bf (ii)}Let $f$ and $g$ be in $C^q(R)$, 
$1\leq q\leq n$, and suppose that 
$$(f, g)_q\overset\text{def.}\to= \sum_{\sigma\in {\frak I}_q^*} f(\sigma )g(\sigma )$$
makes sense, suppose moreover that $q\geq 1$ and that $f=d^{q-1}(h)$, with\break 
$h\in C^{q-1}(R)$, suppose also that
$(h, \delta^q(g))_{q-1}$ makes also sense. Then, it is not difficult to prove that
$$(d^{q-1}(h), g)_q=(h, \delta^q(g))_{q-1}$$ 
then, we see that our next definitions of
laplacians and harmonic cocycles are classical (see \cite{GarH}). 

\proclaim{(3.2.4) Definition} Let $q$ be an integer, $1\leq q\leq n$, the $q$-th laplacian is
the map $\Delta^q\overset\text{def}\to= 
\delta^{q+1}\circ d^q+d^{q-1}\circ \delta^q: C^q(R)@>>>C^q(R)$. We also set 
$\Delta ^0=\delta ^1\circ d^0$. 
\endproclaim

\proclaim{(3.2.5) Definition} Let $q$ be an integer, $0\leq q\leq n$, we say that  an
element $f$ of $C^q(R)$ is an harmonic cocycle of level $q$, or an harmonic
$q$-cocycle, with values in $R$, if $f$ is in the kernel of $\Delta^q$ . 
Let $\underline {H}^q(R)$ be the set of
harmonic cocycles of level $q$ with values in $R$ (then $\underline {H}^q(R)=Ker(\Delta^q)$).\par
Let $\Gamma $ be a subgroup of $G(K_\infty)$. Let also $\underline {H}^q(R)^\Gamma$ be 
the set of $q$-harmonic cocycles which are invariant under the $\Gamma$-action and
$\underline {H}_!^q(R)^\Gamma$ be the set of those which moreover have compact (finite) support
modulo $\Gamma $.
\endproclaim

\proclaim{(3.2.6) Remark} Suppose, as in (3.2.3), that all the expressions written with 
$( \ , \ )_\cdot $ have sense. Take $f\in C^q(R)$, we have
$$(\Delta^q(f), f)_q=(\delta^{q+1}\circ d^q(f), f)_q+(d^{q-1}\circ \delta^q(f), f)_q
=(d^q(f), d^q(f))_{q+1}+(\delta^q(f), \delta^q(f))_{q-1}$$ 
Then, if for instance $R\subset {\Bbb R}$, one has $\Delta^q(f)=0$ if and only if 
$d^q(f)=0$ and $\delta^q(f)=0$.\endproclaim 
  The next proposition give an important case where this characterization is valid.

\proclaim{(3.2.7) Proposition} Let $\Gamma $ be an arithmetic group and let
$R$ be a subring of $\Bbb R$. Let 
$f\in C^q(R)^\Gamma $ ($f$ is stable under $\Gamma $), with $0\leq q\leq n$. 
Suppose that the support of $f$ is finite modulo $\Gamma $, then $f$ is an harmonic cocycle if
and only if $d^q(f)=0$ and $\delta^q(f)=0$.
\endproclaim

\demo{Proof} For any $q$, $1\leq q\leq n$, set $C^q (\Gamma \backslash {\frak I}^\star ,R)
=Hom_{\Bbb Z}({\Bbb Z}[\Gamma \backslash {\frak I}_q^\star], R)$. For any $q$, $1\leq q\leq n$,
 for any $f\in C^q (\Gamma \backslash {\frak I}^\star ,R)$ 
and any $s\in \Gamma \backslash {\frak I}_{q-1}^\star$, let
$$\delta_1^q(f)(s)=\sum_{\sigma \in \Gamma \backslash {\frak I}_q^\star,\sigma >s}
\eta (s,\sigma ){\sharp (\Gamma_s)\over \sharp (\Gamma_\sigma )}f(\sigma )$$
(see (3.2.1)), where $\Gamma_s$ and $\Gamma_\sigma $ are the stabilizers in $\Gamma $ of 
lifts in ${\frak I}$, say $\tilde s$ and  $\tilde \sigma $, of $s$ and $\sigma $, where, if
$\tilde s$ and  $\tilde \sigma $ are chosen such that $\tilde s$ is a face of  $\tilde \sigma $,
one has writen $\eta (s,\sigma )$ in place of $\eta ({\tilde s},{\tilde \sigma })$ . 
We also set $\delta_1^0=0$.\par
 For any $q$, $0\leq q\leq n-1$, for any $f\in C^q (\Gamma \backslash {\frak I}^\star , B)$ and
any $\Sigma \in \Gamma \backslash {\frak I}_{q+1}^\star$ let
$$d_1^q(f)(\Sigma  )=\sum_{\sigma \in \Gamma \backslash {\frak I}_q^\star,\sigma <\Sigma }
\eta (\sigma ,\Sigma )f(\sigma ).$$
We also set $d_1^n=0$.\par
  A function $f\in C^q(R)^\Gamma $ gives rise to one, say $f_1$, in
$C^q (\Gamma \backslash {\frak I}^\star , R)$ (if $\sigma_1$ is the class mod. $\Gamma $ of
$\sigma \in {\frak I}_q^\star$, we have $f_1(\sigma_1)=f(\sigma )$) and it is not difficult to
see that
$$\delta^q(f)=\delta_1^q(f_1) \ \ , \ \ d^q(f)=d_1^q(f_1).$$\par
  For $f,g\in C^q (\Gamma \backslash {\frak I}^\star , R)$, $f$ or $g$ having a finite support, let
$$(f , g)_q^\Gamma =\sum_{\sigma \in \Gamma \backslash {\frak I}_q^\star}
{1\over \sharp (\Gamma_\sigma )}f(\sigma )g(\sigma ).$$ 
Let $1\leq q\leq n$, $f\in C^{q-1} (\Gamma \backslash {\frak I}^\star , R)$ and
$g\in C^q (\Gamma \backslash {\frak I}^\star , R)$, $f$ or $g$ having a finite support, one has
$$(d_1^{q-1}(f),g)_q^\Gamma =(f,\delta_1^q(g))_{q-1}^\Gamma .$$
It follows that, for any $q$, $1\leq q\leq n$, and any $g\in 
C^q (\Gamma \backslash {\frak I}^\star , R)$
$$(\delta_1^{q+1}\circ d_1^q(g) \ + \ d_1^{q-1}\circ \delta_1^q(g) \ , \ g)_q^\Gamma  \
= \ (d_1^q(g),d_1^q(g))_{q+1}^\Gamma  \ + \  (\delta_1^q(g),\delta_1^q(g))_{q-1}^\Gamma  \
. $$
The proof for $q=0$ uses the same arguments.\qed
\enddemo

\subhead 3.3 Harmonic cocycles and cohomology of $\Gamma $\endsubhead

  The aim of this paragraph is to prove the theorem (3.3.1) about the cohomology with
compact supports of an arithmetic group $\Gamma$, which is defined in (3.1.5). Recall that
$p$ is the characteristic of our base field $K$ and that we say that a group has no
$p'$-torsion if any of its torsion elements is of order a power of $p$.

\proclaim{(3.3.1) Theorem} Let $\Gamma $ be an arithmetic group whithout $p'$-torsion; for
any $q$, $0\leq q\leq n$, let $H^q(\Gamma , {\Bbb Z}[1/p])_{\text {tors}}$ be the torsion
subgroup of $H^q(\Gamma , {\Bbb Z}[1/p])$. Then,
we have a canonical isomorphism of ${\Bbb Z}[1/p]$-modules
$$H_!^q(\Gamma , {\Bbb Z}[1/p]) \ \simeq \ {\underline H}_!^q({\Bbb Z}[1/p])^\Gamma 
\oplus H_!^q(\Gamma , {\Bbb Z}[1/p])_{\text {tors}}.$$
Moreover, these modules are finitely generated and  
${\underline H}_!^q({\Bbb Z}[1/p])^\Gamma$ is free (see (3.1.6)).
\endproclaim

{\bf (3.3.2)} Let $R$ be a ring, equals to ${\Bbb Z}[1/p]$ or $\Bbb Q$.
For any $q$, $0\leq q\leq n$, let 
$C_!^q(R)^\Gamma $ be the set of elements of $C^q(R)^\Gamma $ which have finite
supports modulo $\Gamma $. Note that the complex $(C^\cdot (R)^\Gamma  , d^\cdot )$ is
canonically isomorphic to that defined in the proof of (3.2.7).
Let $B_!^q(R)^\Gamma $ and $Z_!^q(R)^\Gamma $ be the intersections 
with $C_!^q(R)^\Gamma $ of respectively $B^q(R)$ and $Z^q(R)$.
Let $f$ be in $C^q(R)^\Gamma $, note that for any simplex
$\sigma$ of $\Gamma \backslash {\frak I}^\star$, the expression $f(\sigma )$ makes sense. 
One has (if $q>1$)
$$\delta^q(f)(s)=\sum_{\sigma \in \Gamma \backslash {\frak I}_q^\star,\sigma >s}
\eta (s,\sigma ){\sharp (\Gamma_s)\over \sharp (\Gamma_\sigma )}f(\sigma ).$$
Let moreover $g$ be in $C^q(R)^\Gamma $ and suppose that one of $f$ or $g$ is of
finite support modulo $\Gamma $, we set (see the proof of (3.2.7))
$$(f , g)_q^\Gamma =\sum_{\sigma \in \Gamma \backslash {\frak I}_q^\star}
{1\over \sharp (\Gamma_\sigma )}f(\sigma )g(\sigma ).$$
Let $1\leq q\leq n$, $f\in C^{q-1}(R)^\Gamma $ and
$g\in C^q(R)^\Gamma $, $f$ or $g$ being of finite support modulo
$\Gamma $; it is easy to see that
$$(d^{q-1}(f),g)_q^\Gamma =(f,\delta^q(g))_{q-1}^\Gamma .$$

\proclaim{(3.3.3) Proposition} For any $q$, $1\leq q\leq n$, let
$B_!^q({\Bbb Z}[1/p])_{\text {``tors"}}^\Gamma$ be the set of\break 
$f\in Z_!^q({\Bbb Z}[1/p])^\Gamma$
such that there exists $a\in \Bbb Z$, $a\not= 0$, with $af\in
B_!^q({\Bbb Z}[1/p])^\Gamma$ (i.e. the inverse image of 
$H^q(\Gamma , {\Bbb Z}[1/p])_{\text {tors}}$ by the canonical morphism\break
$Z_!^q({\Bbb Z}[1/p])^\Gamma@>>>H^q(\Gamma , {\Bbb Z}[1/p])$).
We have the direct sum of ${\Bbb Z}[1/p]$-modules
$$Z_!^q({\Bbb Z}[1/p])^\Gamma = {\underline H}_!^q({\Bbb Z}[1/p])^\Gamma
\oplus B_!^q({\Bbb Z}[1/p])_{\text {``tors"}}^\Gamma.$$
\endproclaim

\demo{Proof} One has ${\underline H}_!^q({\Bbb Z}[1/p])^\Gamma \subset Z_!^q({\Bbb Z}[1/p])^\Gamma$
(see (3.2.7)).

\proclaim {(3.3.4) Lemma} {\bf (i)} ${\underline H}_!^q({\Bbb Z}[1/p])^\Gamma $ is orthogonal to
$B_!^q({\Bbb Z}[1/p])_{\text {``tors"}}^\Gamma $, with respect to the scalar product 
$( \ , \ )_q^\Gamma $ (see (3.3.2));\par
  {\bf (ii)} ${\underline H}_!^q({\Bbb Z}[1/p])^\Gamma $ is a free $\Bbb Z$-module of finite rank. 
\endproclaim

\demo {Proof of the lemma} Let $f\in {\underline H}_!^q({\Bbb Z}[1/p])^\Gamma $, 
$g\in B_!^q({\Bbb Z}[1/p])_{\text {``tors"}}^\Gamma $, $a\in {\Bbb Z}-\{0\}$ such that 
$ag\in B_!^q({\Bbb Z}[1/p])^\Gamma $ and let $h$ be a $q-1$-cochain such that
$d^{q-1}(h)=ag$. One has
$$a(f,g)_q^\Gamma = (f,d^{q-1}(h))_q^\Gamma = (\delta^q(f), h)_{q-1}^\Gamma =0$$
because $f\in Ker(\delta^q)$ (see (3.2.7)). It follows also from (3.2.7) that 
${\underline H}_!^q({\Bbb Z}[1/p])^\Gamma $ is torsion free, which implies (ii).
\qed \enddemo

{\it End of the proof of the proposition (3.3.3).}  
Let ${\Cal G}_q$ be a finite set of generators, over $\Bbb Z$, of $H^q(\Gamma , {\Bbb Z}[1/p])$,
containing a set ${\Cal G}'_q$ of generators of $H^q(\Gamma , {\Bbb Z}[1/p])_{\text {tors}}$ 
and containing a basis ${\Cal B}_q$ of ${\underline H}_!^q({\Bbb Z}[1/p])^\Gamma $. Let 
$\tilde {\Cal G}_q$ be a lift of ${\Cal G}_q$ in $Z_!^q({\Bbb Z}[1/p])^\Gamma $, containing
${\Cal B}_q $. Let $Z$ (resp. $B$) be the subset of elements of
$Z_!^q({\Bbb Z}[1/p])^\Gamma$ (resp. $B_!^q({\Bbb Z}[1/p])_{\text {``tors"}}^\Gamma$) generated by
${\Cal G}_q$ (resp. ${\Cal G}'_q$); set $H={\underline H}_!^q({\Bbb Z}[1/p])^\Gamma$. 
Let $Z_{\Bbb Q}$, $B_{\Bbb Q}$ and $H_{\Bbb Q}$ be the same sets with $\Bbb Q$ in place of
${\Bbb Z}[1/p]$.\par
  Note that the $\Bbb Z$-modules $Z$ and $B$ are free of finite ranks and that $H$ is a submodule 
of $Z$, thus, with the previous lemma, we see that the formula\break
$Z=H\oplus B$ implies the proposition. Indeed, because $H_{\Bbb Q}\cap Z=H$, it remains to prove that
$Z_{\Bbb Q}=H_{\Bbb Q}\oplus B_{\Bbb Q}$.\par
 Let $ B_{\Bbb Q}^\bot$ the subspace of $Z_{\Bbb Q}$ orthogonal to $ B_{\Bbb Q}$ with respect to
the pairing $( \ , \ )_q^\Gamma$. We have $Z_{\Bbb Q}=B_{\Bbb Q}^\bot \oplus B_{\Bbb Q}$.
Now we prove that $H_{\Bbb Q}=B_{\Bbb Q}^\bot$. Let $f\in Z_{\Bbb Q}$ and let
$f=u+d^{q-1}g$ be its decomposition  with respect to the previous direct sum ($u$ is in 
$B_{\Bbb Q}^\bot$ and $g$ is a $(q-1)$-cochain). Let $h$ be a $(q-1)$-cochain:
$$(\delta^q(f-d^{q-1}g), h)_{q-1}^\Gamma = (f-d^{q-1}g, d^{q-1}h)_q^\Gamma 
=(u, d^{q-1}h)_q^\Gamma =0,$$
then $f-d^{q-1}g\in Ker(\delta^q)$ and hence $f-d^{q-1}g\in H_{\Bbb Q}$ (see (3.2.7)). Then
$H_{\Bbb Q}\supset B_{\Bbb Q}^\bot$, the other inclusion is given by (3.3.4).
\qed \enddemo

\noindent {\bf Proof of (3.3.1).} For any $q$, $1\leq q\leq n$,
propositions (3.1.4) and (3.3.3) imply the formula of the theorem, the case $q=0$ is easy.
The other sentence comes from (3.1.6).

\proclaim{(3.3.5) Corollary} Let $\Gamma $ be an arithmetic group whithout $p'$-torsion; for
any $q$, $0\leq q\leq n$, we have a canonical isomorphism of ${\Bbb Q}$-vector spaces
$$H_!^q(\Gamma , {\Bbb Q}) \ \simeq \ {\underline H}_!^q({\Bbb Q})^\Gamma .$$
and these spaces are finite dimensional over $\Bbb Q$.
\endproclaim

\head 4 The locus of supports of harmonic cocycles and of cohomology. \endhead

 Let $\Gamma $ be an arithmetic group whithout $p'$-torsion. The ${\Bbb Z}[1/p]$-modules
${\underline H}_!^\cdot ({\Bbb Z}[1/p])^\Gamma $ and $H_{!}^\cdot (\Gamma ,{\Bbb Z}[1/p])$ are
finitely generated (see (3.3.1)), 
then there exists a finite subcomplex of $\Gamma \backslash {\frak I}$
which contains all the supports modulo $\Gamma $ of elements 
of these modules (the support of an
element of these cohomology modules is defined by (3.1.4) and (3.1.5)). 
The purpose of this paragraph
is to give a description of the locus of this subcomplex. Indeed, we first introduce 
a notion of cusps of an arithmetic group $\Gamma $, they are ``the missing part of 
$\Gamma \backslash {\frak I}$ at infinity", and, secondly, we prove that 
``neighborhoods of this missing part at infinity" are given modulo $\Gamma $ by a union
of finitely many sectors of $\frak I$, and that the complementary complexe of each of these
neighborhoods is finite modulo $\Gamma $. The supports of the elements of 
${\underline H}_!^\cdot ({\Bbb Z}[1/p])^\Gamma $ and of $H_{!}^\cdot (\Gamma ,{\Bbb Z}[1/p])$ are
modulo $\Gamma $ outside sufficiently small neighborhoods of the cusps. These properties are
known in the case where $n=1$ (i.e. when the Building is a tree, see \cite {Se 2} ch. II, \S
2.3).

\subhead 4.1. The cusps\endsubhead

In all this paragraph, $\Gamma $ is an arithmetic subgroup of $G(K)$. Recall that 
${\frak I}_\infty$ is the spherical building at infinity of $\frak I$ (see \S 2.3).

\proclaim{(4.1.1) Definition}  Let $v_0$ be a fixed vertex of $\frak I$. For any chamber
$\sigma $ of $\frak{I}_\infty$, let $\Cal {S}_\sigma $ be a sector of $\frak I$
beginning at $v_0$ and ending at $\sigma $ (see (2.3.1) and (2.3.3)). 
Let $\sigma $ and $\sigma '$ be two chambers of
$\frak{I}_\infty$, we say that they are equivalent modulo $\Gamma $ if
$$sup_{\gamma \in\Gamma } \ \sharp(({\Cal S}_\sigma \cap 
{\Cal S}_{\gamma (\sigma ')})_n)=\infty$$ 
i.e. if the number of chambers of ${\Cal S}_\sigma \cap 
{\Cal S}_{\gamma (\sigma ')}$ is unbounded, as $\gamma $ goes through $\Gamma $.\par
  A $\Gamma$-equivalent class of chambers of  $\frak{I}_\infty$ will be called a cusp of
$\Gamma$. 
\endproclaim

This notion of cusp  does not depend on the choice of the vertex $v_0$. This is an equivalent
relation on the set of chambers of $\frak{I}_\infty$.\par
  
\proclaim{(4.1.2) Proposition} The set  of cusps of $\Gamma $ is finite, it 
canonically one to one with $\Gamma \backslash G(K)/ P(K)$. \endproclaim

\demo{Proof} The second part of the proposition implies the first, because it is wellknown that
\proclaim{(4.1.3) Lemma} The set $\Gamma \backslash G(K)/ P(K)$ is finite. \endproclaim
\noindent (see for instance [Go]).\par
 Let $\sigma_1$,...,$\sigma_d$ be chambers of ${\frak I}_{\infty}$ which
represent
$\Gamma \backslash G(K)/ P(K)$ (see (2.3.6)), 
let ${\Cal S_1}$,...,${\Cal S_d}$ be sectors of $\frak I$ 
which end respectively at $\sigma_1$,...,$\sigma_d$ 
and which begin at the same vertex $v_0$ of $\frak I$.
Let $\Cal S$ be a sector of $\frak I$, begining in $v_0$, and, for any $l>0$, 
let ${\Cal S}(l)$ be the
set of chambers of $\Cal S$ with distance at most $l$ from $v_0$. For all $l$ there exists $i$
and $\gamma $ , there exists a sector ${\Cal S}_{i,\gamma }$ beginning in $v_0$ and ending at
$\gamma(\sigma_i)$, such that ${\Cal S}(l) \subset {\Cal S}_{i,\gamma }$; the conclusion now is
a direct consequence of the definition (4.1.1).\qed
\enddemo

\proclaim{(4.1.4) Remark} Indeed, the preceeding proof means that the cusps of $\Gamma $ are
``rationnal" over the global field $K$ (this is known in dimension $1$: \cite{Se} ch. II, \S
2.3).\par 
  In the sequel, we will say equivalently that a cusp of $\Gamma $ is 
a class of chamber of $\frak{I}_\infty$ or an element of $\Gamma \backslash G(K)/ P(K)$.
\endproclaim

\subhead 4.2 On the structure  of $\Gamma\backslash {\frak I}$\endsubhead

The following result gives the main property of our notion of cusp.
First we introduce notations which will be used in the rest of this chapter.\par

\vskip5pt \noindent {\bf (4.2.1)} let $\Gamma $ be an arithmetic subgroup of $G(K)$.  
Let $v_0$ be a fixed vertex of
the building $\frak I$. Let $\sigma_1$,...,$\sigma_d$ be chambers of ${\frak I}_{\infty}$ 
which represent $\Gamma \backslash G(K)/ P(K)$, i.e. which represent the cusps of 
$\Gamma $.  Let ${\Cal S_1}$,...,${\Cal S_d}$ be sectors of $\frak I$ which
end respectively at $\sigma_1$,...,$\sigma_d$ and which begin at $v_0$. For any integer
$l\geq 0$ and for any $i$, $1\leq i \leq d$, 
let ${\Cal S_i^l}$ be the subcomplex of ${\Cal S_i}$ such that for any chamber $C$ of
${\Cal S_i^l}$, the (combinatorial) distance from $C$ to the complementary complex ${\Cal S_i^c}$ of
${\Cal S_i}$ (in $\frak I$) 
is at least $l$, i.e. such that for any chamber $C'$ of ${\Cal S_i^c}$, any geodesic
from $C$ to $C'$ contains at least $l$ chambers. Let
$$P(l)=\cup_{1\leq i \leq d}\cup_{\gamma \in \Gamma}\gamma({\Cal S_i^l})$$
and $D(l)$ be its complementary complex in $\frak I$.\par
$P(l)$ is, modulo $\Gamma $, what one has called, in the introduction of this paragraph, a 
``neighborhood of this missing part at infinity", this missing part being modulo $\Gamma $ the
union of the cusps.  
Note that $\Gamma $ acts on $P(l)$ and $D(l)$.

\proclaim{(4.2.2) Theorem} For all $l\geq 0$, the complex $\Gamma \backslash D(l)$ is
finite.
\endproclaim

 This theorem implies that $P(l)$ does not depend, up to a finite 
 subcomplex modulo $\Gamma $, on 
the choices  of the $\sigma_i$'s and of the choice
of the vertex $v_0$. One has also the following remark.

\vskip5pt \noindent {\bf (4.2.3)} Suppose that $\Gamma $ is whithout $p'$-torsion, then, for $l$
sufficiently large, the supports of the elements of
${\underline H}_!^\cdot ({\Bbb Z}[1/p])^\Gamma $ and $H_{!}^\cdot (\Gamma ,{\Bbb Z}[1/p])$
are in $D(l)$.\par

\vskip5pt
 All the rest of this \S  4 is devoted to the proof of the theorem (4.2.2).\par 
{\it First we reduce to the case} $\Gamma = G(A)$. Let $\sigma $ be a chamber at infinity which
represents a cusp of $G(A)$ and let ${\Cal S}$ be a sector of $\frak I$ ending at $\sigma $. Let 
$g_1(\sigma ),\cdots ,g_r(\sigma )$, with $g_1,\cdots ,g_r\in G(A)$, be representatives of the
cusps of $\Gamma $, which are equal modulo $G(A)$ to the cusp $\sigma $ of $G(A)$. Let
${\Cal L}_1,\cdots , {\Cal L}_r$ be sectors in $\frak I$ ending respectively at
$g_1(\sigma ),\cdots ,g_r(\sigma )$. There exists a choice of these 
last sectors such that,
for any integer $l\geq 0$, the two complexes
$$\bigcup_{\gamma \in \Gamma, \ 1\leq i\leq r} \gamma({\Cal L}_i^l)
\ \ \  {\text {and}} \ \ \  
\bigcup_{\gamma \in \Gamma, \  1\leq i\leq r} \gamma g_i({\Cal S}^l) $$
(${\Cal L}_i^l$ and ${\Cal S}^l$ are defined as in (4.2.1))
differ modulo $\Gamma $ by finitely many chambers. This implies our assertion, because the number
of cusps is finite.\par

\vskip 5pt \noindent {\bf (4.2.4)}
Then {\it we suppose now that $\Gamma = G(A)$}. We denote by $C_0$ and $v_0$ 
respectively the fundamental chamber and the
fundamental vertex (see (2.2.4)).\par

\vskip 5pt \noindent {\bf (4.2.5)} Let $Y$ be a subset of $G({\Bbb A}_f)$ which represents
the set of double classes $G(K)\backslash G({\Bbb A}_f)/G({\Bbb O}_f)$. It is wellknown that
$Y$ is finite (indeed one to one with Pic$A$). For any ${\underline y}\in Y$ one sets
$$\Gamma_{\underline y} \ = \ {\underline y}G({\Bbb O}_f){\underline y}^{-1}\cap G(K).$$
These groups are arithmetic, we can suppose that 
${\underline 1}=(1,\cdots ,1)\in G({\Bbb A}_f)$ is in $Y$, then we have
$\Gamma_{\underline 1}=G(A)$.\par
It is also wellknown
that there exists a one to one map
$$G(K)\backslash G(\Bbb A)/(G({\Bbb O}_f)\times B_\infty) \rightarrow
\coprod_{{\underline y}\in Y} \Gamma_{\underline y}\backslash {\frak I}_n$$
(where ${\frak I}_n$ is the set of chambers of the building ${\frak I}$) induces by
 $$G({\Bbb A}) = G({\Bbb A}_f)\times G(K_\infty)\rightarrow 
\coprod_{{\underline y}\in Y}\Gamma_{\underline y}\backslash {\frak I}_n $$
which maps ${\underline g}=(\gamma {\underline y}{\underline k}_f,g_\infty)$ to
 $\gamma^{-1}g_\infty \ {\text {mod.}} \ \Gamma_{\underline y}$ (where 
$\gamma \in G(K)$, ${\underline y}\in Y$ and ${\underline k}_f\in G({\Bbb O}_f)$;
${\underline g}$ being given, ${\underline y}$ is well defined).\par

\vskip5pt This last remark explains why some arguments of the proof of
theorem (4.2.2) will be of adelic nature. This is the object of the next
section.

\subhead 4.3 Proof of theorem (4.2.2): the adelic part \endsubhead             

\vskip 5pt \noindent {\bf (4.3.1)} Let $\frak V$ be the set of normalized  valuations of
$K$.    If ${\underline t}=(t_\omega )_{\omega \in \frak V}$ is an id\`ele of $K$, i.e. is an
element of 
$GL_1(\Bbb A)$, we set
$$|{\underline t}|=\prod_{\omega \in \frak V} |t_\omega |_\omega .$$ 
Let $H$ be a subset of $G(\Bbb A)$, let
${\underline g}=({\underline g}_{i, j})$ be an element of $H$; 
for any integer $i$, $0\leq i\leq n$, we set 
$|t_i({\underline g})|=|{\underline g}_{i, i}|$, if it makes sense,
which is supposed in the following definitions.\par
 Let $c$ and $c'$ be two positive real numbers with $c'\leq c$, we set 
$$H_{c'}=\{ {\underline g}\in H \ /\ \forall i \ 0\leq i\leq n-1\
c'\leq |t_i({\underline g})/t_{i+1}({\underline g})|\},$$
$$H_{c'}^c=\{ {\underline g}\in H \ /\ \forall i \ 0\leq i\leq n-1\
c'\leq |t_i({\underline g})/t_{i+1}({\underline g})|\leq c\},$$
$$H_{(c',c)}=\{ {\underline g}\in H \ /\ \forall i \ 
c'\leq |t_i({\underline g})/t_{i+1}({\underline g})|\ , \exists i\ 
|t_i({\underline g})/t_{i+1}({\underline g})|>c \ , \ 0\leq i\leq n-1\}.$$

  Let $\omega \in \frak V$. The preceeding definitions make sense for a subset 
$H$ of $G(K_\omega )$,
if one considers, as usual, $G(K_\omega )$ as a subgroup of $G({\Bbb A})$. In that case, to be
precise, one will write ${}^\omega H_{c'}$, ${}^\omega H_{c'}^c$ and ${}^\omega H_{(c',c)}$ (the
reason of that will be clear when one will consider subsets of $G(K)$, $G(K)$ being
viewed as a subgroup of $G(K_\omega )$).

\proclaim{(4.3.2) Lemma}  
There exits two real numbers $ 0<c'_0<c_0$ such that, for all 
$c'$ and $c$ satisfying $c'\leq c'_0<c_0\leq c$, the image in 
$G(K)\backslash G(\Bbb A)/G({\Bbb O})$ \par
\noindent 1. of $P({\Bbb A})_{c'}$ is onto,\par
\noindent 2. of $P({\Bbb A})_{c'}^c$ is finite.\par
One will suppose moreover (for technical reasons) that $c'\leq 1$.\endproclaim

\demo{Proof} This lemma is a direct consequence of deep results of G. Harder (apply
$g\mapsto g^{-1}$ to S\"atze 2.1.1 and 2.2.2 of [Ha 1], see also \S 1.1 of \cite {Ha 2}).
\qed \enddemo
 
The next proposition is the key result of this section. Before to set it, we introduce a
notation.

\vskip5pt \noindent {\bf (4.3.3)} For any ${\underline y}\in Y$ (see (4.2.5)), let
${\frak S}_{\underline y}$ be a subset of $G(K)$ one to one with the set of double classes
$\Gamma_{\underline y}\backslash G(K)/P(K)$, i.e. a subset of $G(K)$ which represents the
cusps of the arithmetic group $\Gamma_{\underline y}$ (see (4.1.2)).  

\proclaim{(4.3.4) Proposition} 
There exists a constant $\delta >0$ satisfying the following property:
let $c'$ and $c$ be such that $c'\leq c'_0<c_0\leq c$ and $\delta c'\leq c'_0<c_0\leq \delta c$
(see (4.3.2)), then there exists a finite subset $F\subset U(K)$ such that
$$G(K)(P({\Bbb A})_{(c' ,c)})(G({\Bbb O}_f)\times B_\infty )
\subseteq G(K)\bigg [ \coprod_{{\underline y}\in Y} 
(\{ {\underline y}\} \times E_{{\underline y}}) \bigg ]  (G({\Bbb O}_f)\times B_\infty),$$
where for any ${\underline y}\in Y$,  
$E_{{\underline y}}={\frak S}_{\underline y}
F({}^\infty T(K)_{(\delta c', \delta c)}){\overline W}$ (this is a subset of $G(K_\infty)$). 
\endproclaim 

  The constant $\delta $ does not depend on $c'$ and $c$, but $F$ depends on them. The proof
of this proposition needs many lemmata.

\proclaim{(4.3.5) Lemma} $P({\Bbb A}_f)=U(K)T({\Bbb A}_f)P({\Bbb O}_f)$.
\endproclaim

\demo{Proof} Let ${\underline g}\in P({\Bbb A}_f)$, ${\underline g}=({\underline g}_{i,j})$.
Let ${\underline t}=diag({\underline g}_{0,0},\cdots ,{\underline g}_{n,n})\in T({\Bbb A}_f)$.
We want to find $u=(u_{i,j})\in U(K)$ such that 
${\underline w}:={\underline t}^{-1}u{\underline g}\in P({\Bbb O}_f)$. This gives the relations
for $i\leq j$
$${\underline w}_{i,j}=\sum_{1\leq k\leq j}{\underline g}_{i,i}^{-1}{\underline g}_{k,j}u_{i,k}
\in P({\Bbb O}_f).$$
$u_{i,i},\cdots, u_{i,j}$ having been chosen, because $K$ is dense in ${\Bbb A}_f$, we see easely
that $u_{i,j+1}$ exists.
\qed \enddemo

\proclaim{(4.3.6) Lemma} Let $X$ be a set of representatives in $GL_1({\Bbb A}_f)$ of\break
$GL_1(K)\backslash GL_1({\Bbb A})_f/GL_1({\Bbb O}_f)$. There exists a constant $\delta_1>0$
satisfying the following property: let $c$ and $c'$ be as in (4.3.2), 
for any ${\underline g}\in P({\Bbb A})_{(c',c)}$ there exists
${\underline \tau }\in T(X)$ and $g_\infty \in {}^\infty P(K_\infty )_{(\delta_1 c',\delta_1 c)}$
such that
$${\underline g}\equiv ({\underline \tau }, g_\infty) \ mod. \ 
(P(K) \ , \ P({\Bbb O}_f)\times \{1\})$$
(this formula means $mod. \ P(K)$ on the left and 
$mod. \ P({\Bbb O}_f)\times \{1\}$ on the right).
\endproclaim

\demo{Proof} We write, following (4.3.5), 
${\underline g}=(u{\underline \tau }{\underline \rho} ,h_\infty)$
with $u\in U(K)$, ${\underline \tau }\in T({\Bbb A}_f)$, ${\underline \rho} \in P({\Bbb O}_f)$ and 
$h_\infty \in  P({\Bbb O}_\infty)$. We also write  
${\underline \tau }=\gamma {\underline \tau }_1 {\underline \rho }_1$ with $\gamma \in T(K)$,
${\underline \tau}_1\in T(X)$ (see the definition of $X$) and  
${\underline \rho }_1\in T({\Bbb O}_f)$. We have
$${\underline g}\equiv (\tau_1, \gamma^{-1} u^{-1}h_\infty) \ 
mod. \ (P(K) \ , \ P({\Bbb O}_f)\times \{1\})$$
and for all $i=0,\cdots ,n$ (see (4.3.1))
$$\mid t_i({\underline g})\mid =\mid t_i(\tau_1, \gamma^{-1} u^{-1}h_\infty)\mid ,$$
then for all $i=0,\cdots ,n-1$
$$\mid t_i(1, \gamma^{-1} u^{-1}h_\infty)/t_{i+1}(1, \gamma^{-1} u^{-1}h_\infty)\mid =$$
$$\mid t_i(\tau_1, \gamma^{-1} u^{-1}h_\infty)/t_{i+1}(\tau_1, \gamma^{-1} u^{-1}h_\infty)\mid 
.\mid t_i(\tau_1, 1)/t_{i+1}(\tau_1, 1)\mid^{-1}=$$
$$\mid t_i({\underline g})/t_{i+1}({\underline g})\mid . 
\mid t_i(\tau_1, 1)/t_{i+1}(\tau_1, 1)\mid^{-1}$$
We set $\delta_1 =min_{0\leq i\leq n-1, \tau_1\in T(X)}
\mid t_i(\tau_1, 1)/t_{i+1}(\tau_1, 1)\mid^{-1}$, which exists because $X$ is finite. 
\qed \enddemo

\proclaim{(4.3.7) Lemma} There exists a constant $\delta >0$ satisfying the following
property: let $c$ and $c'$ be as in (4.3.2), then for any ${\underline g}\in P({\Bbb
A})_{(c',c)}$ there exists ${\underline y}\in Y$, $\sigma \in {\frak S}_{\underline y}$ 
(see (4.3.3)) and $g_\infty \in ({}^\infty P(K_\infty )_{(\delta c',\delta c)})$ such that
$${\underline g}\equiv ({\underline y}, \sigma g_\infty) \ mod. \ (G(K) \ , \ 
G({\Bbb O}_f)\times \{1\}).$$
\endproclaim

\demo{Proof} We have, following (4.3.6)
$${\underline g}\equiv ({\underline \tau }, g_\infty) \ mod. \ 
(P(K) \ , \ P({\Bbb O}_f)\times \{1\})$$
with, in particular, $g_\infty \in {}^\infty P(K_\infty )_{(\delta_1 c',\delta_1 c)}$. 
There exists ${\underline y}\in Y$, $\gamma \in G(K)$ and ${\underline k}\in G({\Bbb O}_f)$ such
that ${\underline \tau }=\gamma {\underline y}{\underline k}$. Note that $\gamma $ runs in a
finite set depending only on $X$, because $\tau$ is in $T(X)$.\par
  We set $\gamma^{-1}=\alpha \sigma e$ with $\alpha \in \Gamma_{\underline y}$, 
$\sigma \in {\frak S}_{\underline y}$ and $e\in P(K)$ (see (4.3.3)). One has
$${\underline g}\equiv (\gamma {\underline y}{\underline k}, g_\infty) \ mod. \ 
(P(K) \ , \ P({\Bbb O}_f)\times \{1\})$$
$${\underline g}\equiv ({\underline y}, \gamma^{-1}g_\infty) \ mod. \ 
(G(K) \ , \ P({\Bbb O}_f)\times \{1\})$$
$${\underline g}\equiv (\alpha^{-1}{\underline y}, \sigma eg_\infty) \ mod. \ 
(G(K) \ , \ P({\Bbb O}_f)\times \{1\})$$
which implies, because $\alpha \in \Gamma_{\underline y}\subset 
{\underline y}G({\Bbb O}_f){\underline y}^{-1}$,
$${\underline g}\equiv ({\underline y}, \sigma eg_\infty) \ mod. \ 
(G(K) \ , \ G({\Bbb O}_f)\times \{1\}).$$
As $e$ runs in a finite set depending only on $X$ (as the above $\gamma $'s) and $Y$, there exists
$\delta_2>0$, which does not depend on $c'$ and $c$, such that
$eg_\infty \in ({}^\infty P(K_\infty )_{(\delta_1\delta_2 c',\delta_1\delta_2 c)})$.
\qed \enddemo

\proclaim{(4.3.8) Lemma} Let $c$ and $c'$ be as in (4.3.2). Then there
exists a finite subset 
$F$ of $U(K)$ satisfying the following property: for any
${\underline y}\in Y$ and any $\sigma \in {\frak S}_{\underline y}$ (see (4.3.3))
$$\Gamma_{\underline y}\sigma F({}^\infty T(K)_{(c', c)})P({\Bbb O}_\infty) \
\supseteq \ \sigma ({}^\infty P(K_\infty)_{(c', c)}).$$ 
\endproclaim

\demo{Proof} We prove a little better: let $\Gamma $ be an arithmetic group, then there exists a
finite subset $F_\Gamma$ of $U(K)$ such that
$$ \leqno  \text{{\bf (4.3.9)}} \ \ \ ({}^\infty P(K_\infty)_{(c', c)})=(\Gamma \cap
U(K))F_\Gamma  ({}^\infty T(K)_{(c', c)}) P({\Bbb O}_\infty).$$ 
This implies the lemma, for $\Gamma =\sigma^{-1}\Gamma_{\underline y}\sigma $ and because 
$\sigma$ and ${\underline y}$ run in finite sets.\par
  Note that it is sufficient to prove (4.3.9) for $\Gamma =G(A)$, because 
$\Gamma \cap U(K)$ and $G(A)\cap U(K)$ are commensurable.\par

  Let $R\subset K$ be a set of representatives of ${\Bbb F}(\infty)$ (the residue field at
$\infty$), with $0\in R$, let $\eta_0$ be the integer defined by
$$-\eta_0 \ = \ {\text {max}}\{v_\infty (a) \ / \ a\in A, v_\infty (a)<0\},$$
let $a_0\in A$ which realises $-\eta_0$ and set 
($\pi _\infty =\pi \in K$ is uniformising parameter at $\infty$)
$$\tilde R \ = \ \{\rho \pi^{-\eta}a_0^l \ / \  \rho \in R, 0\leq \eta <\eta_0, l\geq 0,
\eta +\eta_0 l>0\},$$
$$ \Pi \ = \ \{\pi^{-\eta}a_0^l \ / \  0\leq \eta <\eta_0, l\in{\Bbb Z}\}.$$

  \par It is easy to see that
$$\leqno \hskip40pt  ({}^\infty P(K_\infty)_{(c', c)})=
({}^\infty T(K)_{(c', c)})U(K_\infty)T({\Bbb O}_\infty),$$
indeed, with straightforward calculations one can prove that
$$\leqno \text{{\bf (4.3.10)}} \ \ \ ({}^\infty P(K_\infty)_{(c', c)})=({}^\infty T(\Pi
)_{(c', c)}) U(\tilde R) P({\Bbb O}_\infty).$$ 
\par
  Let $t={\text {diag}}(t_0,\cdots ,t_n)\in ({}^\infty T(\Pi )_{(c', c)})$ and 
$u\in U(\tilde R)$ and set $tut^{-1}=(u_{i,j})$. 
We have for $i\leq j$ and when $u_{i,j}\not= 0$ 
$$ \mid u_{i,j}\mid_\infty \geq  \mid t_i/t_j\mid_\infty \geq (c')^{j-i}\geq (c')^n.$$ 
On the other hand, we know that there exists a real
number $\varepsilon $ satisfying the following property: for any $\lambda$ in $K$ there exists an
$a$ in $A$ such that
$\mid \lambda -a\mid_\infty \leq \varepsilon $. These two remarks lead to prove that for any
$t\in ({}^\infty T(\Pi )_{(c', c)})$, for any $u\in U(\tilde R)$, there exists 
$\gamma \in U(A)$ such that, if $(v_{i,j})$ denotes the matrix $\gamma^{-1}tut^{-1}$, if
$i<j$ and $v_{i,j}\not=0$, then 
$$\varepsilon \geq \mid v_{i,j}\mid_\infty \geq (c')^n.$$
The set

$$F_1:=\{(v_{i,j}) \ / \ (v_{i,j})=\gamma tut^{-1}, \gamma \in U(A),
t\in ({}^\infty T(\Pi )_{(c',c)}), u\in U(\tilde R),
\varepsilon \geq \mid v_{i,j}\mid_\infty \geq (c')^n \forall i,j \}$$
is finite modulo $P({\Bbb O}_\infty)$ on the right (note that $\tilde R$ and $\Pi $ are discrete);
let $F$ be a finite subset of $U(K)$ such that $F_1\subset FP({\Bbb O}_\infty)$.\par
  Let $g\in ({}^\infty P(K_\infty)_{(c', c)})$ and write following (4.3.10) $g=tuh$, 
with $t\in ({}^\infty T(\Pi )_{(c',c)})$, $u\in U(\tilde R)$ and $h\in P({\Bbb O}_\infty)$;
we have
$$g=tuh=(tut^{-1})tuh=\gamma^{-1}(\gamma tut^{-1})h$$ with $\gamma tut^{-1}\in FP({\Bbb
O}_\infty)$. This proves (4.3.9).
\qed \enddemo

\proclaim{(4.3.11) Lemma} There exists a constant $\delta >0$ satisfying the following
property: let $c'$ and $c$ be such that $c'\leq c'_0<c_0\leq c$ and $\delta c'\leq
c'_0<c_0\leq \delta c$ (see (4.3.2)), then there exists a finite subset $F\subset U(K)$ such
that
$$P({\Bbb A})_{(c',c)}\subseteq  G(K)\coprod_{{\underline y}\in Y}
\bigg [\{{\underline y}\}\times \bigg ({\frak S}_{\underline y}
F ({}^\infty T(K)_{(\delta c',\delta c)})\bigg )\bigg ]
(G({\Bbb O}_f)\times P({\Bbb O}_\infty))$$ 
\endproclaim

\demo { Proof} Let $\delta $ be given by (4.3.7), let
$c$ and $c'$ be such that $c'\leq c'_0$, $\delta c'\leq c'_0$ and $c\geq c_0$, $\delta c\geq c_0$
(see (4.3.2)). Let 
${\underline g} \in P({\Bbb A})_{(c', c)}$. We write as in (4.3.7)
$${\underline g}\equiv ({\underline y}, \sigma g_\infty) \ mod. \ (G(K) \ , \ 
G({\Bbb O}_f)\times \{1\})$$ 
with $g_\infty \in ({}^\infty (P(K_\infty)_{(\delta c', \delta c)})$. \par
  We have, following (4.3.8), $\sigma g_\infty = \gamma \sigma e t \rho$ where 
$t\in ({}^\infty T(K)_{\delta c', \delta c)})$, $e\in F$, $\rho \in P({\Bbb O}_\infty)$ and 
$\gamma \in \Gamma_ {\underline y}$, $\gamma ={\underline y}{\underline k}{\underline y}^{-1}$
with ${\underline k}\in G({\Bbb O}_f)$. Thus we see that modulo
$(G(K) \ , \ G({\Bbb O}_f)\times P({\Bbb O}_\infty))$
$${\underline g}\equiv ({\underline y}, \gamma \sigma e t\rho )
\equiv (\gamma ^{-1}{\underline y}, \sigma e t\rho )=
({\underline y}{\underline k}, \sigma e t\rho )\equiv  ({\underline y}, \sigma e t).\qed $$ 
\enddemo

\vskip5pt \noindent {\bf Proof of proposition (4.3.4).} There exists finitely many
elements $p_i$ of $P({\Bbb O}_\infty)$, $1\leq i\leq i_0$, such that
$$G({\Bbb O}_\infty)=\bigcup_{w\in {\overline W}, 1\leq i\leq i_0}p_iwB_\infty$$
(because the existence of the canonical morphism $G({\Bbb O}_\infty)/B_\infty \hookrightarrow 
G({\Bbb F}(\infty))/P({\Bbb F}(\infty))$),
then we have (as before, ${\underline 1}$ means $(1,\cdots ,1)\in G({\Bbb A}_f)$) 
$$P({\Bbb A})_{(c',c)}\bigg (\{{\underline 1}\}\times G({\Bbb O}_\infty)\bigg )
= P({\Bbb A})_{(c',c)}\bigg (\{{\underline 1}\}\times ({\overline W}B_\infty)\bigg )$$
which gives with (4.3.11) the expected formula.

\subhead 4.4 End of the proof of theorem (4.2.2) \endsubhead 

  To finish the proof, we need also the two following easy lemmata.

\proclaim{(4.4.1) Lemma} Let $l>0$ be an integer and let ${\Cal T}_l$ be the 
set of elements of $T(K)$ which give on $V_0$ (the fundamental appartement)
translations by vectors of the form $\sum_{0\leq i\leq n-1}m_ie_i$ with
$m_i\geq l$ for all $i$ (see (2.2.1)).
\roster
\item Let ${\Cal S}_1$ be the sector of $V_0$ (see (4.2.1) and 2.3) 
beginning at $v_0$ (the fundamental vertex) and ending at the fundamental chamber at
infinity, then
$$({\Cal S}_1^l)_n \subseteq
{\Cal T}_l{\overline W}C_0:=\{twC_0/t\in {\Cal T}_l, w\in {\overline W}\}
\subseteq ({\Cal S}_1^{l-1})_n$$
where ${\Cal S}_1^l$ is defined in (4.2.1), $({\Cal S}_1^l)_n$ 
is the set of chambers of ${\Cal S}_1^l$, $C_0$ is the fundamental chamber, 
and $\overline W$ the linear part of the affine Weyl group $W$.
\item Let $F$ a finite subset of $U(K)$ (as that defined in (4.3.4)). Recall that 
${\frak S}_{\underline 1}\subset G(K)$ 
is a set of representatives of the cusps of  $\Gamma_{\underline 1} =G(A)$, 
i.e. of $G(A)\backslash G(K)/P(K)$; we can suppose that $1$ is in
${\frak S}_{\underline 1}$ and represents the fundamental chamber at 
infinity. Then there exists an integer $\kappa \geq 0$ such that, for 
any $l\geq \kappa$,
$${\frak S}_{\underline 1}F{\Cal T}_l{\overline W}C_0
\subseteq (\cup_{1\leq i\leq d}{\Cal S}_i^{l-\kappa })_n.$$

\endroster
\endproclaim

\demo{Proof} Part (1) is a direct consequence of calculations which prove that
the closed chamber $\overline C_0$ is
a fundamental domain for the action of $G(K_\infty)$ on $\frak I$ 
(see \cite {Bro} ch. I, \S 5.F anf the definition of $\frak I$, as in (2.2.2)).\par
  Part (2) follows from part (1) and from the fact that $U(K)$ stabilizes the fundamental chamber at infinity.
\qed \enddemo

\proclaim{(4.4.2) Lemma} Let $l>0$ be an integer and let $c_1'\leq 1$, $c_1$ be real numbers
such that
$$log(c_1)\geq l(n+1)log(\sharp {\Bbb F}(\infty))-{(n-1)(n+1)^2 \over 4}log(c_1').$$
Then, $({}^\infty T(K)_{(c_1', c_1)})\subset {\Cal T}_l$.
\endproclaim

\noindent {\bf Remark.} There exists such $c_1'\leq 1$ and $c_1$ which moreover satisfy
(4.3.4).\par 

\demo{Proof} Let $t=diag(z_0\pi^{k_0},\cdots ,z_n\pi^{k_n})$ be an element of 
$T(K)_{(c_1', c_1)}$, with $z_i$ in ${\Bbb O}_\infty^\star \cap K$ and $k_i\in {\Bbb Z}$. This
gives on
$V_0$ the translation by the vector $\sum_{0\leq i\leq n-1}m_ie_i$ with, for all $i$,
$m_i=-(k_0+\cdots +k_i)$ (recall that $k_0+\cdots +k_n=0$). 
Because $t$ is in $T(K)_{(c_1', c_1)}$, one has, for all $i$,
$k_i-k_{i+1}\leq -(log(c_1'))/log(\sharp {\Bbb F}(\infty))$ and there exists $i_0$ such that
$k_{i_0}-k_{i_0+1}\leq -(log(c_1))/log(\sharp {\Bbb F}(\infty))$.\par
   Let $i$ be an integer, $0\leq i\leq n-1$. For any integer $j$ with 
$0\leq j\leq i$ let $\alpha _{j, i}=(j+1)(n-i)/(n+1)$ and, when $i<j\leq n-1$, let
$\alpha _{j, i}=(i+1)(n-j)/(n+1)$. Then we have 
$k_0+\cdots +k_i=\sum_{0\leq j\leq n-1}\alpha _{j, i}(k_j-k_{j+1})$.\par
  It follows from all these remarks that for all $i$
$$\align 
m_i  &=-\sum_{0\leq j\leq n-1}\alpha _{j, i}(k_j-k_{j+1})  \\
{}   &\geq \frac1{log(\sharp {\Bbb F}(\infty))}\bigg ( (n-1)(log(c_1'))max_j(a_{j, i})+(log(c_1))
min_j(\alpha _{j, i})\bigg )\\
\endalign$$
because $c_1'\leq 1$ and $c_1\geq 1$; this gives
$$m_i\geq \frac1{(n+1)log(\sharp {\Bbb F}(\infty))}\bigg 
( {(n-1)(n+1)^2 \over 4}log(c_1')+log(c_1) \bigg ) \geq l.\qed $$
\enddemo

  {\it Now we can finish the proof of the theorem}. Let $\delta >0$ 
be the constant given by (4.3.4), 
let $c'$ and $c$ be such that $c'\leq c'_0<c_0\leq c$ and 
$c'_1=\delta c'\leq c'_0<c_0\leq c_1=\delta c$ (see (4.3.2)). 
 Let $l>0$ be an integer. We suppose moreover that $l$, $c'_1$ and
$c_1$ satisfy part (2) of (4.4.1) and (4.4.2). Let $\varphi $ be the canonical map
$$\varphi :{\frak I}_n\rightarrow \Gamma_{\underline 1}\backslash {\frak I}_n$$
(recall that $\Gamma_{\underline 1}=G(A)$).\par
  It follows from (4.4.1) that 
$$\varphi ({\frak S}_{\underline 1}F{\Cal T}_l{\overline W}C_0) \subseteq
\varphi  (P(l-\kappa )_n).$$
 One has also (see the choice of $l$)
$$({}^\infty T(K)_{(c_1', c_1)})\subset {\Cal T}_l.$$
Thus, it is suffiscient to prove that the complementary complex  of\break 
$\varphi ({\frak S}_{\underline 1}F({}^\infty T(K)_{(c_1', c_1)}){\overline W}C_0)$, 
in $\Gamma_{\underline 1}\backslash {\frak I}_n$, is finite. Let
$$G(\Bbb A)@>\eta >> 
G(K)\backslash G(\Bbb A)/(G({\Bbb O}_f)\times B_\infty) @>\psi >>
\coprod_{{\underline y}\in Y} \Gamma_{\underline y}\backslash {\frak I}_n$$ be 
the two maps of (4.2.5), the second, $\psi $, being one to one.
It follows from (4.3.4) that 
$\psi^{-1}\circ\varphi ({\frak S}_{\underline 1}F({}^\infty T(K)_{(c_1', c_1)}){\overline W}C_0)$
is in the complementary set of $\eta (P({\Bbb A})_{(c' ,c)})$ and this last set is finite
(see (4.3.2)). 
\qed

\head 5 Euler-Poincar\'e characteristic.\endhead

  Let $\Gamma $ be an arithmetic group with no $p'$-torsion.
  In this paragraph we prove that a good choice of the locus of supports of elements of
${\underline H}_!^\cdot ({\Bbb Q})^\Gamma $ leads to determine ``geometrically" the
Euler-Poincar\'e characteristic of $\Gamma $ for the cohomology with compact supports and with
coefficients in $\Bbb Q$ (indeed our proof goes for all subfield of $\Bbb R$).
 This result is a  generalization to dimension $n>1$ of some
aspects of \cite {Se 2} ch. I \S 3.3, \cite {Re}, \cite {Ge-Re} \S 3.2 (see \S (1.1.3)). 

\vskip5pt
\noindent {\bf (5.0.1)} Let $\Gamma $ be an arithmetic group. One denotes by 
$\chi _!(\Gamma ,{\Bbb Q})$ its Euler-Poincar\'e characteristic for the cohomology 
with values in $\Bbb Q$ and with compact supports, i.e.
$$\chi _!(\Gamma ,{\Bbb Q})=\sum_{0\leq q\leq n}(-1)^q{\text {dim}}_{\Bbb Q}
H_!^q(\Gamma ,{\Bbb Q}).$$

\vskip5pt \noindent {\bf (5.0.2)} 
Let $l$ be an integer such that any cocycle in
$C_!^\cdot ({\Bbb Q})^\Gamma $ has, up to a coboundary, its support in $D(l)$ (see
(3.3.5), (4.2.2) and (4.2.3)).\par
Let $D\subset {\frak I}$ such that\par
\noindent {\bf (i)} as subspace of $\frak I$, $D$ is contractible and
$\Gamma .D\supseteq D(l)$; \par
\noindent {\bf (ii)} $D$ is a finite subcomplex of $\frak I$,
(to be a subcomplex means that $D$ contains all the faces of its simplices);\par
\noindent {\bf (iii)} any cocycle in
$C_!^\cdot ({\Bbb Q})^\Gamma $ has, up to a coboundary, its support in $\Gamma .D$ and the support
of any element of ${\underline H}_!^\cdot ({\Bbb Q})^\Gamma $ is in $\Gamma .D$.\par

\vskip5pt The aim of this paragraph is to prove the following

\proclaim{(5.0.3) Theorem} Let $\Gamma $ be an arithmetic group without $p'$-torsion and
let $D$ be as in (5.0.2). For any $q$, $0\leq q\leq n$, let 
$D_{q,\Gamma }$ be the set of non oriented $q$-simplices $\sigma $ of $D$, satisfying the
following property: there exists $\gamma \in \Gamma $ such that $\gamma (\sigma )$ is a
$q$-simplex of $D$ and $\gamma (\sigma )\not= \sigma $. Let $g_q=\sharp (D_{q,\Gamma })$. 
Then one has
$$\chi _!(\Gamma ,{\Bbb Q})=1+\sum_{0\leq q\leq n}(-1)^{q+1}g_q.$$
\endproclaim

\demo{Proof} We use notations and definitions of \S 3.
For any $q$, $0\leq q\leq n$, let $C_D^q$, $Z_D^q$, $B_D^q$ be respectively 
the subspace of $C_!^q:=C_!^q(\Bbb Q)^\Gamma $, $Z_!^q:=Z_!^q(\Bbb Q)^\Gamma $, 
$B_!^q:=B_!^q(\Bbb Q)^\Gamma $ (with $B_!^0:=B_!^0(\Bbb Q)^\Gamma :=0$) of elements with
supports in $\Gamma .D$; one sets also $H_!^q:={\underline H}_!^q(\Bbb Q)^\Gamma $
(see (3.3.5)).
We have (see the proof of (3.3.4)) $Z_D^q=H_D^q\oplus B_D^q$ for $0\leq q\leq n$.\par
  Let $q<n$, the map $\delta^{q+1}:C_D^{q+1}@>>>C_D^q$ is defined, its restriction to $B_D^{q+1}$
is injective, indeed, if $d^q(f)$ is in $Ker\delta^{q+1}$ (with $f\in C_!^q$), then (see (3.3.2))
$$0=(\delta^{q+1}d^qf, f)_q^\Gamma =(d^qf , d^qf)_{q+1}^\Gamma $$
thus $d^qf=0$ (here we use that we are in $\Bbb R$). It follows then
$$C_D^q\simeq Z_D^q\times B_D^{q+1}\simeq H_!^q\times B_D^q\times B_D^{q+1}$$
this formula being also right when $q=0, n$ if we set $B_D^0=B_D^{n+1}=0$.\par
  Let $c_q$, $h_q$, $b_q$ be the dimensions over $\Bbb Q$ of respectively
 $C_D^q$, $H_!^q$, $B_D^q$. Let $D_q$ be the set of $q$-simplexes of $D$. We have,
$0\leq q\leq n$, 
$$c_q=\sharp (\Gamma \backslash \Gamma .D_q)=\sharp (D_q)-g_q $$
and, following the previous caculations,
$$\sum_{0\leq q\leq n}(-1)^q c_q = \sum_{0\leq q\leq n}(-1)^q(h_q+b_q+b_{q+1})
= 1+ b_0 + \sum_{0\leq q\leq n}(-1)^q h_q,$$
and $b_0=0$. We have found (see (3.3.5))
$$\chi _!(\Gamma ,{\Bbb Q})=1+\sum_{0\leq q\leq n}(-1)^{q+1}g_q +\bigg (-1 +
\sum_{0\leq q\leq n}(-1)^q\sharp (D_q)\bigg ).$$
The following wellknown lemma finishes the proof.

\proclaim{(5.0.4) Lemma} Let $D$ be a non empty subcomplex of $\frak I$ 
which is, as a subspace, contractible. For any $q$,
$0\leq q\leq n$, let $D_q$ be the set of $q$-simplices of $D$. Then
$$\sum_{0\leq q\leq n}(-1)^q\sharp (D_q) = 1.$$
\endproclaim
\enddemo

\vskip5pt \noindent {\bf (5.0.5)} {\it A problem}.  
Notations are those of theorem (5.0.3). We can suppose that, for any
chamber $C$ of $D$ and for any $\gamma \in \Gamma $, $\gamma (C)$ is not in $D$ (i.e. 
$D_{n, \Gamma }$ is empty). Let $q$ be an integer, $0\leq q\leq n-1$, and $\Gamma _q$ be the subset
of $\gamma \in \Gamma $ satisfying the following property: $\gamma $ is in 
$\Gamma _q$ if there exists a $q$-simplex $\sigma $ of $D$ such that $\gamma (\sigma )$ is again
in $D$. Let $G_0$ be the subgroup of $\Gamma $ generated by $\Gamma _0$ and 
$\Gamma _{\text {tors}}$ ($\Gamma _{\text {tors}}$ is the group generated by the torsion elements
of $\Gamma $) and, 
for any $q$, $1\leq q\leq n-1$, let $G_q$ be the normal subgroup of $G_{q-1}$
generated by $\Gamma _q$ and $\Gamma _{\text {tors}}$. Set $G_n=\{1\}$. Then we can expect
$$ G_0 = \Gamma $$
and
$${\text {Hom}}_{\text {gr}}(G_q/G_{q+1}, {\Bbb Z}[1/p])\simeq 
{\text {Hom}}_{\text {gr}}(\Pi _{q+1}(\Gamma \backslash {\frak I} ,\cdot), {\Bbb Z}[1/p]) \ \ , \ \
0\leq q\leq n-1$$
(we consider homotopy groups for a fixed point of $\Gamma \backslash {\frak I}$). Such formulae
would generalize properties of arithmetic groups in dimension $1$ (see (1.1.3).

\enddocument